\documentclass[11pt]{article}
\usepackage{amsmath,amsthm,amssymb,graphicx}

\newcommand{\HH}{\mathbb H}
\newcommand{\RR}{\mathbb R}

\newcommand{\NN}{\mathbb N}
\newcommand{\QQ}{\mathbb Q}
\newcommand{\CC}{\mathbb C}
\newcommand{\DD}{\mathbb D}

\newcommand{\calA}{\mathcal A}
\newcommand{\calC}{\mathcal C}
\newcommand{\calD}{\mathcal D}
\newcommand{\calF}{\mathcal F}
\newcommand{\calH}{\mathcal H}
\newcommand{\calL}{\mathcal L}
\newcommand{\calO}{\mathcal O}

\newcommand{\calU}{\mathcal U}
\newcommand{\calV}{\mathcal V}

\newcommand{\del}{\partial}
\newcommand{\olM}{\overline{M}}
\newcommand{\IH}{{I\!H}}

\newcommand{\II}{{\mathbb I}{\mathbb I}}

\newcommand{\Fr}{\mathrm{rel}}
\newcommand{\Kr}{\mathrm{abs}}
\newcommand{\mM}{\mathrm{mM}}
\newcommand{\Mm}{\mathrm{Mm}}

\newcommand{\al}{\alpha}

\newcommand{\frakd}{\mathfrak{d}}

\newcommand{\ran}{{\mathrm{ran}}}
\newcommand{\maxdom}{\mathcal{D}(d_{\max},g,a)}

\newtheorem{proposition}{Proposition}

\newtheorem{definition}{Definition}
\newtheorem{theorem}{Theorem}
\newtheorem{lemma}{Lemma}

\newtheorem{claim}{Claim}
\newtheorem{corollary}[theorem]{Corollary}

\newcommand{\frakp}{\mathfrak{p}}
\newcommand{\frakm}{\mathfrak{m}}

\def\slashs{\,\hbox{\slash}\kern-6.5pt S}
\def\slashd{\,\hbox{\slash}\kern-8.0pt D}

\begin{document}

\title{Harmonic forms on manifolds with edges}
\author{Eugenie Hunsicker \thanks{Partially supported by 
the NSF through an ROA supplement to grant DMS-0204730}
\\ Lawrence University \and
Rafe Mazzeo \thanks{Supported by the NSF through the grant
DMS-0204730}\\ Stanford University}

\maketitle   
\begin{abstract}
Let $(X,g)$ be a compact Riemannian stratified space with simple edge 
singularity. Thus a neighbourhood of the singular stratum is a bundle of 
truncated cones over a lower dimensional compact smooth manifold.
We calculate the various polynomially weighted de Rham cohomology spaces
of $X$, as well as the associated spaces of harmonic forms. In the unweighted
case, this is closely related to recent work of Cheeger and Dai \cite{CD}.
Because the metric $g$ is incomplete, this requires a consideration of the various
choices of ideal boundary conditions at the singular set. 
We also calculate the space of $L^2$ harmonic forms for any complete
edge metric on the regular part of $X$.
\end{abstract}

\section{Introduction}
One of the early successes in the extension of Hodge theory to manifolds with singularities 
was the work of Cheeger in the early 1980's on manifolds with isolated  conic singularities 
\cite{Ch1}. This provided the inspiration for, and one of the first corroborations of, 
conjectures made by him and Goresky and Macpherson relating Hodge theory on stratified 
spaces to intersection cohomology. Soon afterwards, Cheeger gave a general strategy to extend 
these results to singular spaces with an iterated stratified structure \cite{Ch2}. 

In the intervening years substantial progress has been made in this subject. The goal of much
of this work is the following: consider a particular class of noncompact or singular spaces, 
and a natural class of metrics on them, and find a relationship between the space of $L^2$ harmonic 
forms and some purely topological invariants of the underlying space. Somewhat nonobviously,
the dimension of this Hodge cohomology space is a quasi-isometry invariant of the metric, 
which makes this a feasible program. The problem is rather different for complete and for 
incomplete metrics since in the incomplete case one must also understand the contributions
coming from the choice of boundary conditions at the singular locus. 

Methods that have been successful for this encompass sheaf-theoretic 
techniques and both soft and hard analytic 
techniques. As an example of the state of the art
of the former we refer to Saper's recent work on locally symmetric spaces \cite{Sa1},
\cite{Sa2}; not surprisingly, combining these various tools can be very effective,
cf.\ our recent work with Hausel \cite{HHM}. We do not attempt, however, to list
the many other significant recent contributions to this area.

The specific problem we study here is to determine the topological meaning of the space of 
weighted $L^2$ harmonic forms on a compact manifold with `incomplete edge' singularities.
More precisely, let $X$ be a compact stratified space with only two strata: an open top-dimensional 
stratum $X^{(n)} := M$ and a stratum $B$ of dimension $b$ in its closure. We assume that
a tubular neighbourhood of $B$ in $X$ is diffeomorphic a bundle of cones over $B$ with fibre
a (truncated) cone $C_1(F)$ over a smooth compact manifold $F$. The Riemannian metric $g$ on 
$M$ restricts to a standard conic metric on each fibre.

We recast this slightly differently. Let $\olM$ be a smooth compact manifold with boundary $Y = \del \olM$. 
Suppose that $Y$ is the total space of a fibration $\phi:Y \to B$ with fiber $F$, and that $x$ is a 
boundary defining function on $\olM$, so $Y = \{x=0\}$ and $dx \neq 0$ there. Write
\[
n= \dim \olM, \quad b=\dim B, \quad f = \dim F.
\]
A metric $g$ on the interior $M = \olM \setminus Y$ is called an incomplete edge metric if in
some collar neighborhood of the boundary $\calU \cong (0,x_0) \times Y$ of $\del M$, 
it is quasi-isometric to one of the form
\[
g = dx^2 + \tilde{h} + x^2 \kappa,
\]
where $\tilde{h}$ is the pullback to $\calU$ of a metric $h$ on $B$ (via the projection $\calU 
\to Y \to B$), and $\kappa$ is a symmetric two-tensor on $\calU$ which restricts to a metric on 
each fiber $F$ in $\del \olM$, i.e.\ at $x=0$. The metric completion of $M$ with respect to such 
a $g$ is diffeomorphic in an appropriate sense to the stratified space $X$ obtained
by collapsing each fibre $F$ at $Y$ to a point.  Natural examples of incomplete edge
metrics include hyperbolic cone metrics with unbranched singular set, cf.\ \cite{HK}. 

Hodge theory on incomplete edge spaces, subject to the condition that $f$ is 
odd, is also the subject of the paper \cite{CD} by Cheeger and Dai. Their aim,
however, is primarily directed at the study of the signature on cone bundles and its
relationship with Dai's $\tau$-invariant for the bundle $Y \to B$. The present paper 
came into being because we realized that the methods developed in our previous work 
\cite{HHM} adapt directly to this setting, and can be used to determine the individual 
Hodge cohomology spaces also when $f$ is even; it requires little extra effort
to study this problem on an entire scale of polynomially weighted spaces. This 
extended setting is quite natural and in \cite{Hu}, the first author ties it to 
another interpretation of Dai's invariant $\tau(Y)$. Of course, we fully 
acknowledge the overlap of the material here with that in \cite{CD}, and are also
grateful to X.\ Dai for several very useful conversations. 

The main results in \cite{HHM} identify the Hodge cohomology on a manifold $M$ with
the same differential topological structure as above, but endowed with a `fibred boundary' 
or `fibred cusp' metric, with certain intersection cohomology groups of the space $X$.
These types of metrics are complete and occur frequently in interesting geometric contexts, e.g.\ as 
gravitational instanton metrics (the fibred boundary case) or 
locally symmetric metrics with $\QQ$-rank $1$ cusp ends (the fibred cusp case).
The proofs there proceed by first showing that the appropriate intersection cohomology 
can be calculated in terms of weighted (conormal) $L^2$ cohomology, and then showing
that these are identified with $L^2$ harmonic forms. The first step uses sheaf-theory
and the second relies on analysis via a parametrix construction. 

While the proofs here are similar, there are two important differences: first, the 
parametrix method in \cite{HHM} draws on the fibred boundary pseudodifferential calculus 
from \cite{MaMe-phi}, while here we use the edge pseudodifferential calculus from \cite{Ma-edge}.
{}From the reader's point of view, this substitution is only formal, since the results
we need appear quite similar (although the intricacies of the parametrix constructions
and analytic phenomena in the two calculi are quite different). However, since the metrics
in one of the classes
we consider in this paper are incomplete, we must pay more careful attention to the whole question
of choices of closed extension for $d$, $\delta$, $D = d+\delta$ and $\Delta$. 

We now describe our results in more detail. We first consider the weighted 
de Rham complex $(x^a L^2\Omega^*(M,g),d)$, where $g$ is an incomplete edge metric.
In general, if $\calF$ is any function space on the Riemannian manifold $(M,g)$, then 
we denote by $\calF\Omega^*(M)$ the space of sections of the exterior bundle 
$\Lambda^*T^*(M)$ with this regularity. When there is dependence 
on the metric, it is indicated explicitly by writing $\calF\Omega^*(M,g)$. 
A form $\alpha$ is in $x^a L^2\Omega^*(M,g)$ if $\alpha = x^a \alpha'$ where
$\alpha' \in L^2\Omega^*(M,g)$. Some results are presented in 
\S 3.3 concerning when $d$ has a unique closed extension to these weighted spaces.
In particular we prove the

\medskip

\noindent{\bf Proposition:}\ (\rm \S 3.3, Corollary 3) 
{\it Suppose that $((f-1)/2-a,(f+1)/2-a) \cap \NN = \emptyset$ or else, if there exists an integer 
$q_a \in ((f-1)/2-a,(f+1)/2-a)$ then $H^{q_a}(F) = \{0\}$. (Note that this is true when $a=0$ and $b$ 
is even, or else when $b$ is odd and $H^{f/2}(F) = \{0\}$.) Under either of these conditions, 
the operator $d$ on $x^aL^2\Omega^*(M,g)$ has a unique closed extension in all degrees.}

\medskip

This result is an analogue of, and generalizes, a result due to Cheeger 
in the conic case \cite{Ch2}. The proof involves the justification of a delicate
integration by parts. 

Although we state and prove this proposition separately, it is also a consequence
of another result we prove later in this paper concerning closed extensions 
of the elliptic operator $D_a = d + \delta_a$, where $\delta_a$ is the codifferential on 
$x^a L^2\Omega^*(M,g)$. The more technical proof in this case uses the ellipticity of
$D_a$ and the existence of a parametrix for it in the edge calculus, as described 
in \S 4. We prove the 

\medskip

\noindent{\bf Proposition:}\ (\rm  \S 4.3, Theorem 7): {\it Let $(M,g)$ satisfy the
hypotheses in the preceding proposition relative to the weight $a$. If in addition $\Delta_F$
(the Laplacian on the fibres $F$ with respect to any one of the family of metrics $\kappa$)
has no `small eigenvalues', as defined in \S 4.3, then $D_a$ is essentially self-adjoint on 
$x^a L^2\Omega^*(M,g)$. If the first hypothesis is satisfied, then it is always possible to 
achieve this extra small eigenvalue hypothesis with a metric $g'$ which is quasi-isometric to $g$.}  

\medskip

These results concern special situations where there is a unique closed extension,
but in general, there are two different canonical procedures to extend $(x^a L^2\Omega^*(M,g),d)$ 
to a Hilbert complex, known as the maximal or minimal extension of $d$, and these may be lead to
quite different complexes. The associated de Rham cohomologies 
are denoted $H^k_{{\max}/{\min}}(M,g,a)$, respectively. In analogy with familiar terminology on manifolds 
with boundary (which is a special case of our setting when the fibre $F$ is a point), we call the Hodge 
Laplacians associated to each of these complexes the absolute and relative Hodge Laplacians, and denote 
their nullspaces by $\calH^*_{\Kr/\Fr}(M,g,a)$. It is almost a tautology that the maximal and minimal 
weighted de Rham cohomology spaces are identified with the absolute and relative weighted Hodge cohomology
spaces. As for the topological interpretation, we prove the

\medskip

\noindent{\bf Theorem:}\ ({\rm \S 3.5, Theorem 4}){\it Let $(M,g)$ be a manifold with incomplete 
edge metric. The maximal and minimal weighted Hodge cohomology spaces are canonically identified 
with intersection cohomology for the stratified space $(X,B)$ by:
\[
H^k_{\max}(M,g,a) =
\left\{ 
\begin{array}{lll}
I\!H^k_{\overline{\frakm} + \ll a-1 \gg} (X,B) & \qquad & f \mbox{ odd} \\
I\!H^k_{\overline{\frakm}  + \ll a-1/2 \gg } (X,B) & \qquad & f \mbox{ even} 
\end{array}
\right.
\]
and
\[
H^k_{\min}(M,g,a) = 
\left\{ 
\begin{array}{lll}
I\!H^k_{\underline{\frakm} + <a>} (X,B) & \qquad & f \mbox{ odd} \\
I\!H^k_{\underline{\frakm} + <a-1/2>} (X,B) & \qquad & f \mbox{ even} 
\end{array}
\right.;
\]
\noindent
here $\ll t \gg $ denotes the least integer strictly greater than $t$ and $<t>$ denotes the 
least integer greater than or equal to $t$, and $\underline{\frakm}$, respectively $\overline{\frakm}$,
are the lower and upper middle perversities.}

We single out two important special cases:

\medskip

\noindent{\bf Corollary:} {\it The maximal and minimal de Rham cohomologies when $a=0$ 
correspond to upper and lower middle perversity intersection cohomology.
\begin{equation}
\begin{array}{rcl}
H^k_{\max}(M,g,0) &= &  I\!H^k_{\overline{\frakm}}(X) \\
H^k_{\min}(M,g,0) & = &  I\!H^k_{\underline{\frakm}}(X)
\end{array}
\end{equation}
Moreover, when $f$ is even, the maximal and minimal de Rham
cohomologies at weights $\pm 1/2$ coincide, and again correspond
to upper and lower middle perversity intersection cohomology.
\begin{equation}
\begin{array}{rcl}
H_{{\max}/{\min}}^k(M,g,-1/2) &=& I\!H^k_{\overline{\frakm}}(X)\\
H_{{\max}/{\min}}^k(M,g,1/2)  & = & I\!H^k_{\underline{\frakm}}(X).
\end{array}
\end{equation}}

\medskip

The notation $\IH^*_{\frakp}(X,B)$ is somewhat nonstandard, and indicates a slight 
generalization of these spaces (so as to include, for example, the case where $F$ 
is a point), which we discuss in \S 3.2. 

These results are, to some extent, `soft' in that they do not require any
serious use of elliptic theory, and for that reason we have separated
them into the first few sections of the paper. 
The main ingredients in their proofs are some abstract 
functional analytic results involving Hilbert complexes from
\cite{BL}, reviewed in \S 2, the sheaf-theoretic characterization of intersection 
cohomology from \cite{CGM}, discussed in \S 3.2, and the 
appropriate Poincar\'e Lemmas (also known as `the local calculations'), which
are developed in \S 3.4 and 3.5. Very helpful in our approach is the fact that we 
may restrict attention to conormal forms, but this is not strictly speaking
necessary. This part of the paper is a recapitulation and extension of 
Cheeger's original work on the Hodge theory on cones, with an attempt to present
the argument as cleanly as possible in this slightly more general context, but the results
could all also be proved using the techniques in \cite{Ch2}.

On the other hand, we require more analytic information in order to study the
minimal Hodge cohomology $\calH^*_{\min}(M,g,a)$, which is defined as the common 
nullspace of $d_{\min,a}$ and $\delta_{\min,a}$. This terminology is slightly
unfortunate, since $\calH^*_{{\min}}$ is not related to the minimal de Rham
cohomology $H^*_{{\min}}$ discussed above, but rather corresponds to the
nullspace of the minimal extension of $D_a$. In any case, we prove the 

\medskip

\noindent{\bf Theorem:}\ (\rm \S 4.5, Theorem 8) {\it Let $M$ be a manifold with an incomplete
edge metric, $g$.  The minimal weighted Hodge cohomology is given by
\[
\calH^{k}_{{\min}}(M,g,a) =
\left\{ 
\begin{array}{rclll}
\mbox{Im}\,\left(I\!H^k_{\underline{\frakm} + <a>}(X,B) \right.
& \to &  \left. I\!H^k_{\overline{\frakm} + \ll a-1 \gg }(X,B)\right) & \quad & f 
\mbox{odd} \\
\mbox{Im}\,\left(I\!H^k_{\underline{\frakm} + <a-1/2>} (X,B) \right. & \to  & \left.
I\!H^k_{\overline{\frakm} + \ll a-1/2 \gg } (X,B)\right) & \quad & f \mbox{even}.
\end{array}
\right.
\]
In particular, when $a=0$, 
\[
\calH^{k}_{{\min}}(M,g,0) = \mbox{Im\,}(\IH^{k}_{\underline{\frakm}}(X,B) 
\longrightarrow \IH^{k}_{\overline{\frakm}}(X,B)).
\]}

\medskip

The proof requires two main analytic results: the conormal regularity of solutions 
in the minimal domain of the equation $D_a \omega = 0$, and the solvability (and
regularity theory for the solution) of $D_a \zeta = \eta$ for suitable $\eta$.
For these we invoke the theory of pseudodifferential edge operators, as 
developed in \cite{Ma-edge}. The results from this theory which we require 
are reviewed in \S 4.1. 

Refering to that section for the following terminology, we note that the specific
computations we must make in order to apply this more general theory are the calculation 
of the indicial roots of $D_a$ and the injectivity of the normal operator $N(xD_a)$ 
on suitable weighted $L^2$ spaces. As we show in \S 4.2.1, the calculation of indicial 
roots for $D_a$ on manifolds with edges is essentially identical to that on cones, 
and this partially explains the similarity of the results in the two cases. However, 
the extra role played by the model operator $N(xD_a)$, see \S 4.2.2 as well as 
Proposition 8 in \S 4.1, is not required in the conic case, but is the key fact needed
in the parametrix construction in the edge calculus. 

We now turn to some applications and extensions of our results. The first is a Bochner-type
vanishing result. Recall that the Weitzenb\"ock formula for the Hodge Laplacian
on $k$-forms on $M$ states that $\Delta_k = \nabla^* \nabla + {\mathcal R}_k$,
where ${\mathcal R}_k$ is a curvature operator, acting by endomorphisms on
$\bigwedge^kM$. 

\medskip

\noindent {\bf Theorem:} {\it Suppose that $M$ admits an incomplete edge 
metric $g$ such that ${\mathcal R}_k \geq 0$ everywhere, and is strictly positive
at some point of $M$. Then the minimal Hodge cohomology $\calH^k_{{\min}}(M,g)$
(at weight $a=0$), and hence 
\[
\mbox{Im\,}(\IH^{k}_{\underline{\frakm}}(X,B) 
\longrightarrow \IH^{k}_{\overline{\frakm}}(X,B)),
\]
both vanish. If in addition $f$ is odd or else if $f$ is even but
$H^{f/2}(F) = \{0\}$, then we may also deduce that 
$\IH^k_{\underline{\frakm}}(X,B) = \IH^k_{\overline{\frakm}}(X,B) = \{0\}$.}

\medskip

The proof is the usual one, and simply involves noting that when $\omega$ is
in the nullspace of $D_{{\min}}$, then the integration by parts
\[
\langle \Delta \omega, \omega \rangle = ||\nabla \omega||^2 + 
\langle {\mathcal R}_k \omega, \omega \rangle
\]
is justified. 

We can extend this type of analysis significantly further. For example, 
the generalization to second order natural geometric operators $L$ is
essentially straightforward. By definition, such an operator is one 
of the form $L = \nabla^* \nabla + {\mathcal R}$ acting on sections of 
some subbundle $E$ of the full tensor bundle over $M$, using the induced
Levi-Civita connection (though we may also twist by any other bundle with 
connection); the symmetric 
endomorphism ${\mathcal R}$ on $E$ is a generalized curvature operator. 
$L$ is formally symmetric on $L^2(M;E)$, and we may ask the same questions 
about its domain, mapping properties and nullspace as we have for the 
Hodge Laplacian. The point we wish to make is that answers to these 
questions are readily deduced, and direct consequences of the edge theory,
once one has calculated the indicial roots of $L$ and determined the mapping
properties of the normal operator $N(x^2 L)$. Furthermore, the calculation
of indicial roots for such an operator reduces directly to the analogous
computation on the cone $C(F)$ endowed with the metric $dx^2 + x^2
\kappa$. The presence of an extra twisting bundle may change the
arithmetic of the indicial root computation and the spaces on which
the normal operator is injective, but makes no difference elsewise. 
We shall not develop these remarks further here, but shall return to them 
in greater detail elsewhere. (At that time we shall also give a more careful 
account of parametrices in the edge calculus for the Laplacian with relative and 
absolute boundary conditions.)
There are many interesting geometric consequences of such results. 
In particular, the infinitesimal rigidity of hyperbolic cone metrics with smooth 
singular set, as proved in \cite{HK} in $3$ dimensions and in the recent paper 
\cite{Mo} (for deformations amongst Einstein metrics rather than just hyperbolic 
metrics) in dimensions greater than $3$ is a direct consequence.  See \cite{MW}
for further discussion.  

Now consider a more topological application: When $f = \dim F$ is odd, 
the maximal and minimal cohomologies agree, as do the two 
middle perversities $\overline{\frakm}$ and $\underline{\frakm}$ for $X$.  Our results 
in this case, when $a=0$, agree with those in \cite{CD}. In this special case,
and assuming that $n = 4\ell$, there is a well-defined pairing on middle 
degree forms; Cheeger and Dai prove the corresponding signature theorem. We extend 
their result to cover also the case $f$ even. Define the $L^2$-signature as the 
signature of the degenerate pairing on $\calH^{2\ell}_{\Fr}(M,g)$ induced from the map
\[
\calH^{2\ell}_{\Fr}(M,g) \longrightarrow \calH^{2\ell}_{\Kr}(M,g),
\]
and the nondegenerate pairing between these spaces. Similarly, one can
also define a topological signature $\sigma(M)$ as the signature of
the degenerate pairing on $H^*_0(M)$ defined through its map to $H^*(M)$ and 
the nondegenerate pairing between these spaces.  Finally, recall the invariant
$\tau(Y)$ defined by Dai \cite{Dai} for the total space of the fibration $Y = \del M \to B$. 

\medskip

\noindent{\bf Theorem:} {\it  The $L^2$-signature of the stratified Riemannian space $X$
endowed with an incomplete edge metric $g$, is given by
\[
L^2-{\rm sgn}\,(M,g)= \sigma(M) + \tau.
\]}

\medskip

Our final result concerns the Hodge cohomology of the manifold $M$ endowed
with a {\it complete} edge metric. By definition, and following the notation
above, $g$ is a complete edge metric if near $\del M$ it has the form
\[
g = \frac{dx^2 + \tilde{h}}{x^2} + \kappa.
\]
The prototype would be the product of a hyperbolic space (or any conformally compact metric) 
and a compact manifold $F$. It is of interest, as a generalization of the main result of 
\cite{Ma-conf.cpt.}, to calculate the Hodge cohomology for such manifolds. We obtain the

\medskip

\noindent {\bf Theorem:}\ (\rm \S 5, Theorem 9){\it  Let $(M,g)$ be a manifold with a complete edge metric.
Let $X$ be the associated compact stratified space. Suppose that $k$ is 
{\bf not} of the form $j + (b+1)/2$ where $\calH^j(F) \neq \{0\}$. Then
\[
L^2\calH^k(M,g) \cong {\IH}^k_{f+\frac{b}{2}-k}(X,B).
\]
In this case, the $L^2$ signature result is the same as above.
In all other cases, where $k$ does have this form, $L^2\calH^k(M,g)$ is 
infinite dimensional.}

\medskip

In conclusion, let us remark that sorting out the detailed analysis of the Hodge Laplacian 
(and other natural geometric elliptic operators) for more general classes of stratified 
spaces, e.g.\ algebraic varieties, is a problem of great importance with many applications. 
Hodge theory on quite general real analytic manifolds is the subject of an 
ongoing project by D.\ Grieser and R.\ Melrose, using an approach closely related to 
(although more general than) the one used here. We appreciate their
interest in and forebearance concerning the present paper. 

\section{Hodge-de Rham theory for Hilbert complexes}
In this section we review some generalities about $L^2$ cohomology, 
based on the formalism of Hilbert complexes from \cite{BL}, to
which we refer in the interests of brevity for most of the proofs. 

Consider a complex of the form
\begin{equation}
0 \to L_0 \stackrel{D_0}{\longrightarrow} L_1 
\stackrel{D_1}{\longrightarrow}L_2 \ldots \stackrel{D_{n-1}}{\longrightarrow} L_n \to 0,
\label{eq:hc}
\end{equation}
where each $L_i$ is a separable Hilbert space, $D_i:L_i \to L_{i+1}$ is
a closed operator with dense domain $\calD(D_i)$ such that $\ran(D_i) \subseteq 
\calD(D_{i+1})$ and $D_{i+1}\circ D_i = 0$ for all $i$. Under these
conditions, (\ref{eq:hc}) is called a Hilbert complex, and is denoted by
$(L_*,D_*)$.  

Many familiar constructions in Hodge-de Rham theory carry over
immediately to this setting, and we list in particular:
\begin{itemize}
\item[i)] There is a dual Hilbert complex 
\[
0 \to L_0 \stackrel{D_0^*}{\longleftarrow} L_1 
\stackrel{D_1^*}{\longleftarrow}L_2 \ldots \stackrel{D_{n-1}^*}{\longleftarrow} L_n \to 0
\]
defined using the Hilbert space adjoints of the differentials, $D_i^*: L_{i+1} \to L_i$;
\item[ii)] The Laplacian $\Delta_i = D_{i}^* D_i + D_{i-1}D_{i-1}^*$
is a self-adjoint operator on $L_i$ with domain
\[
\calD(\Delta_i) = \{u \in \calD(D_i) \cap \calD(D_{i-1}^*): D_i u \in \calD(D_i^*),
D_{i-1}^* u \in \calD(D_{i-1})\} 
\]
and nullspace
\[
\ker \Delta_i := \calH^i(L_*,D_*) = \ker D_i \cap \ker D_{i-1}^*;
\]
\item[iii)] There is a weak Kodaira decomposition
\[
L_i = \calH^i \oplus \overline{\ran D_{i-1}} \oplus \overline{\ran{D_i^*}};
\]
\item[iv)] The cohomology of $(L_*,D_*)$ is defined by
\[
H^i(L_*,D_*) = \ker D_i/\ran D_{i-1};
\]
if this space is finite dimensional, then $\ran D_{i-1}$ is necessarily
closed and \[
H^i(L_*,D_*) = \calH^i(L_*,D_*).
\] 
\end{itemize} 

The main case of interest here is when $(M,g)$ is a (not necessarily complete) 
Riemannian manifold, $L_i = L^2\Omega^i(M,g)$ and $D_i$ is the exterior derivative 
operator. Later we shall also consider the somewhat more general case 
\[
L_i = e^{2w_i}L^2\Omega^i(M,g),
\]
where $w_i \in \calC^\infty(M)$ is some weight function, but for the
remainder of this section, to be concrete, we shall assume that $w_i \equiv 0$;
it will be clear that everything extends to the general case 
in a straightforward manner.  

To turn the `core' de Rham complex $(\calC^\infty_0\Omega^*(M),d)$ into a Hilbert complex, 
we must specify a closed extension of $d$, and there may be more than one way to do this.

\begin{definition} The two canonical closed extensions of $d$ are:
\begin{itemize}
\item The maximal extension $d_{\max}$; this is the operator 
$d$ acting on the domain
\begin{eqnarray*}
\calD(d_{\max}) &=& \{\omega \in L^2\Omega^*(M): d\omega \in L^2\Omega^*(M,g)\} \\
& = & \{\omega \in L^2\Omega^*(M): \exists \, \eta \in L^2\Omega^*(M,g)   \\
&\ &\quad \mbox{s.t.}\ \langle \omega, \delta \zeta\rangle = \langle \eta, 
\zeta\rangle \ \forall \zeta \in \calC^\infty_0\Omega^*(M)\}.
\end{eqnarray*}
In other words, $\calD(d_{\max})$ is the largest set of forms $\omega$ in $L^2$ 
such that $d\omega$, computed distributionally, is also in $L^2$. 
\item The minimal extension $d_{\min}$; this is given by the graph closure of $d$ 
on $\calC^\infty_0\Omega^*(M)$, i.e.
\begin{eqnarray*}
\calD(d_{\min}) = \{\omega \in L^2\Omega^*(M): \exists\, \omega_j \in
\calC^\infty_0\Omega^*(M),\quad \omega_j \to \omega \ \mbox{in}\ L^2\} \\
\mbox{and}\ d\omega_j \ \mbox{also converges to some}\ \eta \in
L^2\},
\end{eqnarray*}
in which case $d_{\min}\omega = \eta$. 
\end{itemize} 
Maximal and minimal extensions of $\delta$ are defined in the same manner. 
\end{definition}

Clearly $\calD(d_{\min}) \subseteq \calD(d_{\max})$. An old result due to 
Gaffney \cite{dR} shows that these domains are the same when $(M,g)$ 
is complete, but in many other cases of interest (for example, on a 
manifold with boundary) they may differ. 

In order to show that $(L^2\Omega^*(M,g),d_{{{\max}/{\min}}})$ 
are both Hilbert complexes, we require the

\begin{lemma}
\[
d_{\max}: \calD(d_{\max}) \to \calD(d_{\max})
\]
and
\[
d_{\min}: \calD(d_{\min}) \to \calD(d_{\min}).
\]
\label{le:hcdde}
\end{lemma}
\begin{proof} The fact that $(d_{\min})^2 = 0$ follows from the
identity $d^2=0$ on $\calC^\infty_0\Omega^*$, while the fact that 
$\delta^2=0$ on test forms and the definition of distributional 
derivatives shows that $(d_{\max})^2=0$.
\end{proof}
The cohomologies of these complexes are denoted $H^*_{{\max}/{\min}}(M,g)$,
respectively. Here and in the sequel we shall use notation
like ${\max}/{\min}$ in a hopefully self-explanatory manner
to indicate statements which hold for each of the indicated
extensions. 

It is straightforward that the Hilbert complex adjoint of 
$(L^2\Omega^*(M,g),d_{{{\max}/{\min}}})$ is $(L^2\Omega^*(M,g),
\delta_{{{\min}/{\max}}})$, i.e.
\[
(d_{\max})^* = \delta_{\min}\quad \mbox{and}\quad (d_{\min})^* = \delta_{\max}.
\] 

There are three well-behaved weak Kodaira decompositions:
\begin{equation}
L^2\Omega^j(M,g) = \calH^j_{{\Kr}/{\Fr}/{\max}}(M,g) 
\oplus \overline{\ran\, d_{{\max}/{\min}/{\min},\,j-1}} \oplus 
\overline{\ran\, \delta_{{\min}/{\max}/{\min},\, j}}, \label{eq:wkd}
\end{equation}
with summands mutually orthogonal in each case.  The first summand 
on the right, called the absolute, relative or maximal Hodge 
cohomology, respectively, is defined as the orthogonal complement 
of the other two summands. Since $(\ran\, d_{\max})^\perp = \ker\, \delta_{\min}$, 
etc., we see that
\begin{equation}
\calH^j_{{\Kr}/{\Fr}/{\max}}(M,g) = \ker\, d_{{\max}/{\min}/{\max},\, j} \cap
\ker\, \delta_{{\min}/{\max}/{\max},\, j-1},
\label{eq:cint}
\end{equation}
respectively. The third decomposition, incorporating both $d_{\min}$ and $\delta_{\min}$, 
is the original one defined by Kodaira. The corresponding Hodge cohomology $\calH^j_{\max}$ 
is often infinite dimensional, though. We do {\it not} consider a fourth 
weak Kodaira decomposition involving the ranges of $d_{\max}$ and $\delta_{\max}$, since 
these subspaces might not even be disjoint, let alone orthogonal. (This is due to
the fact that $d_{\min}d_{\max}$ may not even be defined, let alone vanish.) 
Nonetheless we still define the minimal Hodge cohomology
\begin{equation}
\begin{array}{rcl}
\calH^j_{\min}(M,g) & = &  L^2 \Omega^j(M,g) \ominus \left( \ran\, d_{{\max},\, j-1}
+  \ran\, \delta_{{\max},\, j} \right) \\
& = & \ker\, d_{{\min},\, j} \cap \ker\, \delta_{{\min},j-1} \\
& = &  \ker\, d_{{\min},\, j}/\, \left(\overline{\ran\, d_{{\max},\, j-1}} \cap
\ker\, d_{{\min},\, j}\right). 
\end{array}\label{eq:hc3}
\end{equation}
Note also that
\[
\calH^j_{\min}(M,g) = \calH^j_{\Fr}(M,g) \cap \calH^j_{\Kr}(M,g)
\]

The operators $d_{{\min}/{\max}}$ are both clearly quasi-isometry 
invariants, and the various Kodaira decompositions above then show that
the minimal and maximal cohomologies $H^*_{{\max}/{\min}}(M,g)$, and 
their `reduced' versions, the absolute and relative Hodge cohomologies 
$\calH^*_{\Kr/\Fr}(M,g)$, are all quasi-isometry invariants. This
invariance is also true for the maximal and minimal Hodge cohomologies
$\calH^*_{{\max}/{\min}}(M,g)$. 
 
There are quite a few `Laplacians' one might consider, most
prominent amongst which are the 
absolute and relative Laplacians
\[
\Delta_{\Kr} = \delta_{\min} d_{\max} + d_{\max} \delta_{\min},\qquad
\Delta_{\Fr} = \delta_{\max} d_{\min} + d_{\min} \delta_{\max}.
\]
These are self-adjoint and satisfy
\begin{equation}
\calH^j_{\Kr}(M,g) = \ker \Delta_{\Kr} \qquad \mbox{and} 
\qquad \calH^j_{\Fr}(M,g) = \ker \Delta_{\Fr}.
\end{equation}
Furthermore, if $H^j_{{\max}/{\min}}(M,g)$ is finite dimensional, then the range of 
$d_{{\max}/{\min},j-1}$ is closed, and $H^j_{{\max}/{\min}}(M,g) = \calH^j_{\Kr/\Fr}(M,g)$. 
Consequently, these Hodge cohomology spaces may be computed using only tools from
differential topology and general cohomology, e.g.\ sheaf theory, Mayer-Vietoris, etc. 

One can also define
\[
\Delta_{\mM/\Mm} = \delta_{{\min}/{\max}} d_{{\max}/{\min}} + 
d_{{\min}/{\max}} \delta_{{\max}/{\min}}.
\]
These are symmetric, but not necessarily self-adjoint, invariant
under the Hodge star, and satisfy
\[
\calH^j_{{\max}/{\min}}(M,g) = \ker\, \Delta_{\mM/\Mm}.
\]
Note that 
\[
\calH^j_{\min}(M,g) = \ker\, d_{\min} \cap \ker\, \delta_{\min},
\]
so this is consistent with our prior definition of $\calH^j_{{\min}}$.

We conclude this section by stating two more results, both true in
the general Hilbert complex setting, but for simplicity we restrict
to the setting of differential forms. The first concerns a K\"unneth-type theorem.

\begin{proposition}[\cite{BL} Corollary 2.15]
Let $(L',D')$ and $(L'' D'') $ be two Hilbert complexes.  Form the 
completed tensor product Hilbert complex $(L,D)$:
\[
L_j = \bigoplus_{i+\ell=j} L'_i \, \hat \otimes \, L''_\ell,
\]
\[
D_j = \bigoplus_{i+\ell=j} (D'_i \otimes \mbox{id}_{L''_\ell} + (-1)^i 
\mbox{id}_{L'_i} \otimes D''_\ell).
\]
\noindent
Suppose that $D''$ has closed range in all degrees. Then
\[
H^j(L,D)= \bigoplus_{i+\ell=j} H^i(L',D') \otimes H^\ell(L''D'').
\]
\label{pr:kunneth}
\end{proposition}

The other result concerns the possibility of computing one of these 
cohomology groups using a `core subcomplex' of smooth (but not 
necessarily compactly supported) forms
\[
\calD^\infty_{{\max}/{\min}} \Omega^*(M,g) \subset L^2 \Omega^*(M,g)
\]
consisting of all elements $\omega$ which are in the domain of $\Delta_{\Kr/\Fr}^\ell$ 
for every $\ell \geq 0$.   
\begin{proposition}[\cite{BL} Theorem 2.12]
The cohomology $H^*_{{\max}/{\min}}(M,g)$ is equal to the cohomology of the
complex $(\calD^\infty \Omega^*_{{\max}/{\min}}(M,g), d_{{\max}/{\min}})$.
\label{pr:core}
\end{proposition}

For example, when $(M,g)$ is compact without boundary, this is simply
the well-known result that $H^*(M)$ can be computed using the complex
of smooth forms. When $(M,g)$ is compact with boundary, then as 
discussed carefully in \cite{BL}, $H^*_{{\max}/{\min}}(M,g)$ is equal 
to the cohomology of the complex of smooth forms continuous to the 
boundary which satisfy absolute/relative boundary conditions. 

\section{De Rham theory and edges}
The context in which we shall adapt and develop the material from
the last section is the category of manifolds with edge
singularities. After defining these we briefly review the 
intersection cohomology theory for such spaces, and then
turn to an analysis of the maximal and minimal cohomologies 
of the Hilbert complexes of weighted $L^2$ forms, and in
particular the identification of different weighted
de Rham cohomologies with intersection cohomologies with
different perversities. 

\subsection{Manifolds with edge singularities}

We now begin to develop some of the ideas in the last section in the
concrete setting of manifolds with conic or edge singularities.

\begin{definition} A pseudomanifold $X$ of dimension $n$ has simple 
edge singularities if it has a dense open stratum $M$, which is
a smooth manifold of top dimension, and the singular strata
$X_{\mathrm{sing}} = X\setminus M$ are a disjoint union of
closed smooth manifolds $B_j$ (of possibly varying dimension)
such that each $B_j$ has a neighbourhood
$\calU_j$ which is diffeomorphic to a bundle with base $B_j$
and fibre a truncated cone $C_1(F_j)$ over a smooth link $F_j$.
\end{definition}

The boundary $\del \calU_j$ of each cone bundle neighbourhood 
$\calU_j$ is the total space of a 
bundle over $B_j$ with fibre $F_j$. There are more complicated
singular spaces with iterated edge singularities, which is why we 
call this class `simple'. However, for brevity, in this paper we shall
refer to a space of this type as a manifold with edge singularities. 
Note that this class includes the case of manifolds with conic singularities,  
i.e.\ where some $B_j$ are $0$-dimensional. 

\begin{definition} 
A metric $g$ on a space $X$ with simple edge singularities is said to
be if incomplete edge type if it is an ordinary smooth metric away
from the singular strata $B_j$, while in each cone bundle neighbourhood
$\calU_j$ it is quasi-isometric to one of the form
\[
g = dx^2 + \pi_j^* h + x^2 \kappa;
\]
here $x$ is the polar distance on each cone $C_1(F_j)$,
$\kappa$ is a symmetric $2$-tensor on $\del \calU_j$ which restricts
to a metric on each fibre $F_j$, $\pi_j: \calU_j \to B_j$ is 
the projection, and $h$ is a metric on $B_j$. 

A metric $g$ on the principal stratum $M$ of such a space $X$ is 
of complete edge type if in each $\calU_j$ it has the form
\[
g = \frac{dx^2 + \pi_j^* h}{x^2} + \kappa,
\]
where $x$, $\kappa$, $\pi_j$ and $h$ are as above. We often drop the $\pi_j^*$ in this 
notation and shall also frequently write $\tilde{g}(x) = h + x^2 \kappa$. 
\end{definition}

We next recall some analytic and geometric properties of Riemannian
submersions from \cite{HHM} and discuss their relevance to de Rham
theory for edge metrics.
 
Let $\phi:Y \to B$ be a fibration with fibre $F$, and suppose that it is endowed
with a metric $\tilde{g}$ of the form $\phi^*(h) + \kappa$, where $h$ is a metric 
on $B$. We assume furthermore that $\phi: (Y,\tilde{g}) \to (B,h)$ is a Riemannian submersion.
The tangent bundle $TY$ splits into a vertical and horizontal subbundle, 
$T^V Y \oplus T^H Y$, where $T^V Y =\mbox{ker\,}(d\phi)$ and $T^H Y$ is its
orthogonal complement (and also the subbundle annihilated by $k$). This induces
a splitting of the form bundles on $Y$, and thus every differential form has a 
(horizontal,vertical) bidegree, i.e.
\[
\Omega^{p,q}(Y) = \Omega^p(B)\, \widehat{\otimes} \, \Omega^q(Y,T^VY).
\]
The space of harmonic forms on $F$ is finite dimensional, and we let
\[
\Pi_0^q: L^2 \Omega^q (F) \longrightarrow L^2 \calH^q(F),
\qquad \Pi_\perp = I - \Pi_0
\]
denote the natural orthogonal projectors; these extend naturally to each 
$L^2\Omega^{p,q}(Y)$. 

\begin{proposition} 
The differential and codifferential on $Y$ decompose as 
\[
d_Y = d_F + \tilde{d}_B - \II + {\rm R}, \qquad 
\delta_Y = \delta_F + (\tilde{d}_B)^* - \II^* + {\rm R}^*,
\]
where $d_F$ is the pullback of $d$ to the fibre, $\tilde{d}_B$ is the 
lift of $d_B$ as a horizontal operator, and $\II$ and ${\rm R}$ are tensorial 
operators built from the second fundamental form of the fibres and the curvature 
of the bundle, respectively. These act as
\[
\begin{array}{rcl}
d_F: \Omega^{p,q}(Y) \to \Omega^{p,q+1}(Y), & \qquad &
\tilde{d}_B: \Omega^{p,q}(Y) \to \Omega^{p+1,q}(Y)\\
\II: \Omega^{p,q}(Y) \to \Omega^{p+1,q}(Y), &\qquad&
{\rm R}: \Omega^{p,q}(Y) \to \Omega^{p+2,q-1}(Y).
\end{array}
\] 
\label{pr:dYdY}
\end{proposition}

Now consider the degenerating family of metrics $\tilde{g}_x = h + x^2\kappa$ 
($0 < x \leq 1$). If $\alpha$ is a $(p,q)$-form, then 
\[
|\alpha|_{\tilde{g}(x)} = x^{-q}|\alpha|_{\tilde{g}(1)}.
\]
Furthermore, as explained in \cite{HHM}, $\tilde{d}_B^x = \tilde{d}_B$, $\II^x = \II$
and $R^x = xR$; the $x$ in the superscript signifies that the operator is to
be calculated relative to the metric $\tilde{g}(x)$, and the operator without
a superscript is calculated relative to $\tilde{g}(1)$. Hence
\[
d_Y^x = d_F + \tilde{d}_B - \II + x{\rm R},\qquad
\delta_Y^x =  \delta_F + (\tilde{d}_B)^* - \II^* + x {\rm R}^*.
\]

We can define the operator
\[
\frakd = \Pi_0\big(\tilde{d}_B - \II\big)\Pi_0;
\]
this acts on the space of fibre-harmonic forms, or
equivalently, we can think of this as acting on the space of forms on $B$
with coefficients in the flat vector bundle of harmonic forms on $F$ that comes
from the fibre bundle $Y$.  In \cite{HHM} we proved the following useful lemmas:
\begin{lemma} The operator $\frakd$ and its adjoint $\frakd^*$
are differentials, i.e.\ $\frakd^2 = (\frakd^*)^2 = 0$.
\label{pr:frakddiff}
\end{lemma}

\begin{corollary} Let $\DD = \frakd + \frakd^*$, and suppose that
$\DD \, \alpha = 0$ for some fibre-harmonic form $\alpha$. Then
$\frakd \alpha = \frakd^* \alpha = 0$, and so the terms $\alpha_{p,q}$
of pure bidegree also satisfy $\DD \alpha_{p,q} = 0$.
\end{corollary}

\subsection{Intersection cohomology}
Let $X$ be a pseudomanifold which is a smoothly stratified space of real dimension $n$, 
with no codimension one stratum.  For this subsection only, we allow $X$ to be 
more singular than was considered earlier in this paper; namely, around any point 
$q\in X$ contained in the stratum $X_{\ell}$ of codimension $\ell$ is a neighbourhood
$\calU_q$ diffeomorphic to $\calV \times C(\Sigma)$, where $\calV$
is a Euclidean ball and $C(\Sigma)$ is the cone over a link $\Sigma$, 
which itself is a stratified space (of dimension $\ell-1$).

A perversity $\frakp$ is an $n$-tuple of natural numbers, $(p(1),p(2), \ldots, p(n))$ 
satisfying $p(1)=p(2)=0$ and $p(\ell-1) \leq p(\ell) \leq p(\ell-1)+1$ for all $\ell \leq n$. 
Associated to such a space $X$ and perversity $\frakp$ is the intersection
complex $I\!C_*^\frakp(X)$; roughly speaking, the integer $p(\ell)$
regulates the dimension of the intersection of chains (in general position)
with the stratum of codimension $\ell$. The homology of this complex
is the intersection homology $I\!H^\frakp_*(X)$; the cohomology of the dual 
cochain complex is the intersection cohomology $I\!H^*_\frakp(X)$.

The following result, which asserts that the cohomology of any fine sheaf 
over $X$ is equal to the intersection cohomology of $X$ (with respect
to some perversity $\frakp$, so long as the local sheaf cohomology in any 
sufficiently small neighbourhood equals the intersection cohomology of that neighbourhood.

\begin{proposition}[\cite{GM2}] Let $X$ be a stratified space
and let $(\calL^*, d)$ be a complex of fine sheaves on $X$ with
cohomology $H^*(X,\calL)$. Suppose that if $\calU$ is a neighbourhood
in the principal (smooth) stratum of $X$, then $H^*(\calU,\calL) =
H^*(\calU,\CC)$, while if $q$ lies in a stratum of codimension
$\ell$, and $\calU = \calV \times C(\Sigma)$ as above, then
\begin{equation}
H^k(\calU,\calL) \cong I\!H_\frakp^k(\calU) =
\left\{ \begin{array}{ll}
I\!H_\frakp^k(\Sigma) & k \leq \ell-2-p(\ell) \\
0         & k > \ell-1-p(\ell).
\end{array} \right.
\label{eq:localcalc}
\end{equation}
Then there is a natural isomorphism between the hypercohomology
$\mathbb{H}^{\,*}(X,\mathcal{L}^*)$ associated to this complex
of sheaves and $I\!H^*_\frakp(X)$, the intersection cohomology of
perversity $\frakp$.
\label{pr:shchar}
\end{proposition}
The details and proof of this theorem can be found in  \cite{CGM} and \cite{Bo}.
We refer to equation (\ref{eq:localcalc}) by saying that $\calL$ satisfies the correct {\it local 
calculation} for intersection cohomology with perversity $\frakp$. This result will
be one of our primary tools below.

On any pseudomanifold $X$, 
there are two distinguished perversity functions:
$\underline{\frakm}$ the lower middle, and $\overline{\frakm}$ the upper
middle perversity. When all strata of $X$ are even dimensional, the two corresponding
intersection cohomologies are the same and satisfy Poincar\'{e} duality. When not
all the strata are all even dimensional, then $I\!H_{\underline{\frakm}}(X) \neq
I\!H_{\overline{\frakm}}(X)$ in general, but these two cohomologies are
Poincar\'{e} dual to one another. 

In the case of an $n$-dimensional manifold with edge singularities, there is 
only one relevant value for the perversity function $\frakp$ in the neighborhood of any of the
components of the singular stratum, namely, its value on the codimension of that stratum,
$X^{sing} \cong B$.
If the fibre $F$ for this stratum has dimension $f$ then this codimension is $f+1$.
If $f$ is odd, then the upper and the lower middle perversities satisfy 
$\underline{\frakm}(f+1) = \overline{\frakm}(f +1) = (f-1)/2$.  If $f$ is even,
then $\underline{\frakm}(f+1) = f/2$ and $\overline{\frakm}(f +1) = f/2 -1$.
We extend the standard definition of perversity and intersection cohomology in this
situation to include perversities with $\frakp \leq 0$ and $\frakp \geq f$.  We
see from the local calculation above that these 
perversities give cohomologies 
$$
\begin{array}{ll}
I\!H^*_\frakp(X, B) \cong H^*(X, X^{sing}) & \frakp(f+1) \geq f \\
I\!H^*_\frakp(X, B) \cong H^*(X- X^{sing}) & \frakp(f+1) \leq 0.
\end{array}
$$
\noindent
When we extend the definition this way, we retain the Poincar\'{e} duality property of 
intersection cohomology:
$$
I\!H^{n/2 - k}_{\underline{\frakm} + s} \cong I\!H^{n/2 + k}_{\overline{\frakm} - s},
$$
\noindent
but we lose the property that the intersection cohomology groups are independent
of stratification, which was the original reason for restricting the values of the perversity
function.  Since there is a natural stratification in the case of manifolds with edge singularities,
this does not cause us trouble here.

\subsection{Maximal and minimal domains for weighted complexes}
{}From now on, we assume that $X$ is a space with simple edge singularities and 
$g$ is a metric on it which is of incomplete edge type. Let us fix a smooth 
function $x$ on the principal stratum $M$ which equals 
the polar distance function in a neighbourhood $\calU_j$ of
each singular stratum $B_j$. 

For any $a \in \RR$, consider the complex of weighted $L^2$ forms
\[
0 \to x^a L^2\Omega^0(M,g) \overset{d}{\longrightarrow} x^a L^2 \Omega^1(M,g) 
\overset{d}{\longrightarrow} \ldots \overset{d}{\longrightarrow}
x^a L^2 \Omega^n(M,g) \to 0.
\]
As already explained in \S 2, there are (at least) two ways to turn this
into a Hilbert complex, using either the maximal or minimal extensions of $d$.
The corresponding complexes, $(x^a L^2\Omega^*(M,g),d_{{\max}/{\min},a})$
are called the maximal and minimal weighted cohomologies (with weight $a$),
respectively, and have cohomologies and spaces of harmonic forms 
$H^*_{{\max}/{\min}}(M,g,a)$ and $\calH^*_{\Kr/\Fr}(M,g,a)$. 
Note that the formal adjoint of $d$ with respect to the 
$x^a L^2$ inner product has a term of order $0$ depending on $a$. 

It will be possible to simplify many of the calculations in the computations
of these spaces below by applying Proposition~\ref{pr:core}. As we show later, 
it will ultimately suffice to consider forms which are `tangentially
regular' or conormal, as we now define. Let $M$ be a compact manifold with boundary. 
The space $\calV_b(M)$ of $b$-vector fields on $M$ is, by definition, the space of 
all smooth vector fields which are constrained at $\del M$ to be tangent to the boundary. 
Thus, in any local coordinate system $(x,y)$ , where $x$ is a boundary defining 
function and $y$ is a local coordinate system on the boundary, $\calV_b$
is generated over $\calC^\infty(M)$ by the basis sections $x\del_x$, $\del_{y_j}$. 

\begin{definition} Let $\gamma \in \RR$. The space of conormal functions on 
$M$ of ($L^\infty$) order $\gamma$, $\calA^\gamma(M)$, is defined as
\[
\left\{u: V_1 \cdots V_\ell \, u \in x^\gamma L^\infty(M) \quad 
\forall \, \ell\ \mbox{and}\ V_j \in \calV_b\right\}.
\]
We write $\calA^{\gamma+}(M) = \cup_{\epsilon > 0} \calA^{\gamma+\epsilon}(M)$
and $\calA^{*}(M) = \cup_{\gamma} \calA^{\gamma}(M)$. 
If $E$ is any smooth vector bundle over $M$, then $\calA^\gamma(M;E)$ is defined using 
any system of smooth trivializations of $E$; in particular, the space of differential 
forms conormal of weight $\gamma$ is denoted $\calA^{\gamma}\Omega^*(M)$.
\end{definition}
Note in particular that if $u \in \calA^\gamma(M;E)$, then $|u| = {\mathcal O}(x^\gamma)$
along with all its $b$ derivatives, and is $\calC^\infty$ in the interior of $M$.
In the following, we frequently make use of the fact that if $u \in 
\calA^*(M;E)$ and $\int |u|^2 x^b\, dxdy < \infty$ (with respect to
any smooth nonsingular Hermitian metric on $E$), then 
$|u| = o(|x|^{-(b+1)/2})$ uniformly in $y \in \del M$ as $x \to 0$.
The proof is elementary and left to the reader. 

We remark that $\calD(d_{{\min},a}) \cap \calA^*\Omega^*$ is dense in 
$\calD(d_{{\min},a})$, and similarly $\calD(d_{{\max},a}) \cap \calA^*\Omega^*$
is dense in $\calD(d_{{\max},a})$. The former of these is immediate from
the definition, while the latter follows (in the cone or edge setting)
using standard mollification arguments. 

Let us now consider the problem of when $\calD(d_{{\max},a})$ equals $\calD(d_{{\min},a})$ 
on the truncated cone $C_1(F)$, with exact warped product conic metric $g = dx^2 + 
x^2 \kappa $. This is already contained in Cheeger \cite{Ch2} (when $a=0$), but we
present this argument to make it accessible for later generalization. By the
remark above, it suffices to consider only conormal forms, and this simplifies 
the discussion somewhat. 

\begin{lemma}  The form $\omega = \alpha + dx \wedge \beta \in 
\calD(d_{\max},a)\cap \calA^*\Omega^k(C_1(F))$ lies in 
$\calD(d_{\min,a}) \cap \calA^*\Omega^k(C_1(F))$ if and only if 
$|\Pi^k_0 \alpha(x)|_g =o(x^{-k})$ (or equivalently, $|\Pi^k_0 \alpha(x)|_\kappa = o(1)$)
when $k \in ((f-1)/2 - a, (f+1)/2 - a)$. In particular, if $k$ is in this range, 
and if $\Phi \in \calH^k(F)$ and $\sigma(x) \in \calA^*(\RR^+)$, then $\sigma(x)\Phi 
\in \calD(d_{\min,a})$ if and only  if $\sigma(x) = o(1)$.
\label{le:max=min}
\end{lemma}

Before commencing with the proof, we state the immediate and important consequence. 
\begin{corollary}
Suppose either that the interval $((f-1)/2-a,(f+1)/2-a)$ contains no
integer value, or else that $k_a \in ((f-1)/2-a,(f+1)/2-a)$ and
$\calH^{k_a}(F) = \{0\}$. Then the extensions $d_{{\max},a}$ and 
$d_{{\min},a}$ coincide.
\end{corollary}

\begin{proof} Following the definitions, to show that the minimal and maximal domains 
are equal it suffices to show that for any conormal forms $\omega \in \calD(d_{{\max},a})$ 
and $\psi \in \calD(\delta_{{\max},a})$,  with $\psi$ vanishing near $x=1$,
we have $\langle d\omega, \psi \rangle_a = \langle \omega, \delta_a \psi \rangle_a$;
in other words, we must show that when $k \in ( (f-1)/2 - a, (f+1)/2 - a))$, 
the boundary term in the integration by parts vanishes if and only if 
$|\Pi^k_0 \alpha(x)|_\kappa = o(1)$.

Define $*_a$ by $\psi \wedge *_a \psi = |\psi|^2 x^{-2a} dV_g$. Then 
$*_{-a}*_a = \pm 1$ (depending on the parity of the degree), 
and we find that $\delta_a = \pm *_{-a} d *_a$ and in addition, if $\psi \in 
x^a L^2\Omega^*(C_1(F),g)$, then $*_a \psi \in x^{-a}L^2\Omega^*(C_1(F),g)$. 

Now
\[
\langle d \omega, \psi \rangle_a = \int d\omega \wedge *_a \psi  = 
\int d(\omega \wedge *_a\psi ) - \int \omega \wedge *_a \delta_a \psi,
\]
and so, writing $*_a \psi = \tau + dx\wedge\rho$, integrating by parts 
produces the boundary term
\begin{equation}
\int_{C_1(F)}  d(\omega \wedge *_a\psi) = \lim_{x \rightarrow 0} 
\int_{\{x\} \times F} \alpha \wedge \tau.
\label{eq:boundary}
\end{equation}
Hence $\omega \in \calD(d_{{\min},a})$ if and only if this boundary term vanishes 
for all $\psi$.

Next, the weighted $L^2$ conditions
\[
\int |\alpha|^2_h \, x^{f-2k-2a} + |\tau|^2_h \, x^{f-2(f-k) + 2a}\, dx \, dV_h < \infty
\]
and conormality yield
\[
|\alpha(x)|_\kappa = o(x^{k - f/2 + a - 1/2}), \qquad |\tau(x)|_\kappa = o(x^{(f-k) - f/2 - a-1/2}),
\]
but this gives only that (\ref{eq:boundary}) is $o(x^{-1})$, which is not enough.

To proceed further, use the Hodge decomposition on $F$
\begin{eqnarray}
\begin{array}{rcl}
\alpha(x) & = & d_F A_1(x) + \delta_F A_2(x) + \alpha_0(x) \\
\tau(x) & = & d_F B_1(x) + \delta_F B_2(x) + \tau_0(x).
\end{array}
\label{eq:hdF}
\end{eqnarray}
Inserting these expressions into (\ref{eq:boundary}), many terms vanish and 
we are left with
\begin{equation}
\lim_{x \rightarrow 0} \left( \int_{F} A_1(x) \wedge \delta_F d_F B_1(x) 
+ \int_{F} d_F \delta_F A_2(x) \wedge B_2(x) + \int_{F} \alpha_0(x) \wedge \tau_0(x)\right).
\label{eq:newboundary}
\end{equation}
Since $A_1, B_1 \perp \ker\, d_F$, we can estimate $||A_1(x)||_{L^2(F)}\leq c 
||\alpha||_{L^2(F)}$, $||B_1(x)||_{L^2(F)} \leq ||\tau||_{L^2(F)}$; in addition,
$d_F \alpha (x) = d_F \delta_F A_2(x) \in x^a L^2 \Omega^{k+1}$,
$d_F \tau (x) = d_F \delta_F B_2(x) \in x^{-a}L^2\Omega^{f-k+1}$.
Hence the first two terms here are $o(1)$, and it remains only to analyze the third.

Now decompose $\alpha_0(x) = \sum a_j(x)\phi_j$ and $*_F \tau_0(x) = \sum t_j(x)\phi_j$, 
where $\{\phi_j\}$ is an orthonormal basis for $\mathcal{H}^k(F)$, and where each $a_j(x)$ 
and $t_j(x)$ is conormal on $[0,1]$. The boundary integral now reduces to the finite sum
$\sum_j a_j(x) t_j(x)$, and we must show that each $a_j(x) t_j(x) = o(1)$. However, 
\[
\int_0^1 \left(a_j^2(x) + (a_j'(x))^2\right) x^{f-2k-2a} \, dx < \infty, 
\qquad 
\int_0^1 \left(t_j^2(x) + (t_j'(x))^2\right) x^{f-2k+2a} \, dx < \infty, 
\]
and these imply that $a_j(x), a_j'(x) = o(x^p)$, $p = k - f/2 + a -1/2$,
and $t_j(x), t_j'(x) = o(x^q)$, $q = (f-k) - f/2 - a -1/2$. The improvement
comes by writing $a_j$, $t_j$ as integrals of $a_j'$, $t_j'$, respectively,
where the other limit of integration is taken at $0$ or $1$ depending on $p$ 
or $q$; this gives
\[
a_j(x) = \left\{
\begin{array}{ll}
o(x^{p+1}) & p < -1 \\
o(\log x) & p = -1 \\
a_j(0) + o(x^{p+1}) & -1 < p < 0 \\
o(x^{p+1})& p \geq 0
\end{array}
\right.,  \qquad
t_j(x) = \left\{
\begin{array}{ll}
o(x^{q+1})& q<-1 \\
o(\log x) & q = -1 \\
t_j(0) + o(x^{q+1}) & q = -1 \\
o(x^{q+1}) & q \geq 0
\end{array}
\right. .
\]
(In particular, $a_j(0)$, $t_j(0)$ exist when $p,q \in (-1,0)$.)
Hence, if either $p \notin [-1,0)$ or $q \notin [-1,0)$, we get $a_j(x)t_j(x) = o(1)$.
If $p=-1$ then $q=0$ and we reach the same conclusion.
If both $-1 < p,q <0$, which occurs precisely when $k \in ( (f-1)/2 - a, (f+1)/2 -a)$,
then we obtain a well-defined, but not necessarily vanishing, boundary term
$\int_F \alpha_0 \wedge \tau_0$.  This completes the proof.
\end{proof}

We wish to extend this result to incomplete edge metrics. Using the quasi-isometry
invariance and localizability (proved in \S 3.4 below) of the maximal and minimal
domains, it suffices to consider forms on $\calV \times C_1(F)$, $\calV \subset
\RR^b$, with warped product metric $dx^2 + h + x^2 k = dx^2 + \tilde{g}(x)$, 
and we may as well assume that $h$ is the Euclidean metric too. As before, $\Pi_0$ 
denotes the orthogonal projection onto fibre-harmonic forms, and we let $Y = \calV \times F$. 
\begin{lemma} Let $\omega = \alpha(x) + dx \wedge \beta(x)$, where
\[
\alpha(x) = \sum_{k} \alpha_{*,k}(x),\quad \beta(x) = \sum_{k} \beta_{*,k}(x)
\]
are the decompositions into fibre degree. Then $\omega \in \calD(d_{{\min},a}) \cap \calA^*$ 
if $\omega \in \calD(d_{{\max},a})\cap \calA^*$ and $|\Pi_0 (\alpha_{*,k}(x))|_{\tilde{g}(x)}
=o(x^{-k})$ whenever $k \in ((f-1)/2-a, (f+1)/2 - a)$. In particular, if $\eta \in 
\Omega^{(*,k)}(Y)$ is fibre harmonic and $s(x)$ is conormal, then $s(x)\eta \in 
\calD(d_{{\min},a})$ if and only if $s(x) = o(1)$.
\label{lem:condom}
\end{lemma}

\begin{proof}
Following the same proof as above, the form $\omega = \alpha + dx \wedge \beta \in \calD(d_{{\max},a}) 
\cap \calA^*$ is also in $\calD(d_{{\min},a})$ if and only if, for every $\gamma \in 
\calD(\delta_{{\max},a})\cap \calA^*$ with $*\gamma = \xi + dx \wedge \nu$, the boundary term 
\[
\int_M d\omega \wedge *\gamma - \int_M \omega \wedge \delta \gamma = 
\lim_{x\rightarrow 0} \int_Y \alpha(x) \wedge \xi(x)
\] 
vanishes. Decomposing into fibre degree, the boundary term becomes:
\[
\lim_{x \rightarrow 0} \sum_{k=0}^f \int_Y \alpha_k(x) \wedge \xi_{f-k}(x).
\]
We have 
\[
\alpha_k \in \calA^{k -(f+1)/2 + a + }(M,\Lambda^{*,k}T^*Y), \qquad
\xi_{f-k} \in \calA^{(f-k) -(f+1)/2 -a + }(M,\Lambda^{*,f-k}T^*Y),
\]
so  $|\alpha_k \wedge \xi_{f-k}|_{\tilde{g}(x)} = o(x^{-1})$ only at this stage.
Continuing as before, we can reduce to the case where each $\alpha_k$
and $\xi_{f-k}$ are fibre-harmonic, and write
\[
(\alpha_{k})_0 = \sum a_j(x,y)\phi_j, \qquad (\xi_{f-k})_0 = \sum t_j(x,y)\phi_j
\]
where $a_j, t_j \in \calA^*$. In order to improve the order of vanishing of
$a_j$ and $t_j$ as before, we observe that, for example, assuming each $\alpha_k$
is fibre-harmonic, then
\[
\omega, d\omega \in x^a L^2 \Rightarrow \int \left(|\alpha_k|^2 + |\beta_k|^2
+ |\del_x \alpha_k - d_y \beta_k|^2 \right)\,
x^{f-2k-2a}\, dx\, dV_{\tilde{g}(1)} < \infty.
\]
However, conormality already gives that 
\[
\int |d_y \beta_k|^2\, x^{f-2k-2a}\, dx\, dV_{\tilde{g}(1)}< \infty,
\]
and so we conclude that
\[
\int |\del_x \alpha_k|^2 \, x^{f-2k-2a}\, dx\, dV_{\tilde{g}(1)} < \infty.
\]
A similar argument applies to $\xi$. The rest of the proof is now the
same as in the conic case.
\end{proof}

\begin{corollary} Suppose that either $((f+1)/2-a,(f-1)/2-a) \cap \NN = \emptyset$
or else if $q_a \in ((f+1)/2-a,(f-1)/2-a) \cap \NN$ then $H^{q_a}(F) = \{0\}$.
(This is true in particular when $a=0$ either when $b$ is even, or else when 
$b$ is odd but $H^{f/2}(F) = \{0\}$.) Then $\calD(d_{{\max},a}) = 
\calD(d_{{\min},a})$, and hence $H^*_{{\max}}(M,g,a) = H^*_{{\min}}(M,g,a) = 
\calH^*(M,g,a)$ in every degree. 
\label{cor:trivext}
\end{corollary}
 
\subsection{Preparations for the Poincar\'{e} lemma}
We now prove several preliminary results which will be used in the 
computation of weighted de Rham cohomology on manifolds with edge 
singularities. Our arguments closely follow those in \cite{Ch2}, but
with simplifications since by Propositions \ref{pr:kunneth} and \ref{pr:core}
it suffices to work only with conormal forms.

To begin, define two complexes of sheaves, $\calL_{{\max}/{\min},a}$;
these are induced by the presheaves obtained by restricting the forms in
$\calD(d_{{\max}/{\min},a})$ to the cofinal sequence of coverings of $X$ of 
the form $\{\calU_\alpha\}$, where each $\calU_\alpha$ is either disjoint from all 
of the singular strata $B_j$ or else, if $\calU_\alpha \cap B_j \neq \emptyset$, then
$\calU_\alpha$ is a product neighbourhood $\calV_{\alpha} \times C_{\epsilon}(F_j)$;
here $\calV_{\alpha} \subset B_j$ and $C_{\epsilon}(F_j)$ is the truncation
to $x \leq \epsilon$ of the cone $C(F_j)$.
\begin{lemma}
The sheaves in each of the two complexes $\calL_{{\max}/{\min},a}$ over $X$ are fine. 
\end{lemma}
\begin{proof} We define a partition of unity $\{\chi_\alpha\}$ subordinate to 
$\{\calU_\al\}$ in such a way that each $d\chi_\alpha$ is bounded. Indeed, away
from the singular strata this is obvious, whereas if $\calU_\alpha = \calV \times 
C_\epsilon(F)$ then we can let $\chi_\alpha$ be a product of cutoff functions in 
each factor. It is now standard that if $\omega \in \calD(d_{{\max}/{\min},a})$, 
then $\chi_\al \omega \in \calD(d_{{\max}/{\min},a})$ as well. This gives the
result.
\end{proof}

Now form the associated spectral sequence for either of the double deRham/\v{C}ech complexes 
$\calL_{{\max}/{\min},a}$. Its hypercohomology may be computed taking either the \v{C}ech 
or the de Rham differential first. In the former case, at the first stage this
becomes the standard definition of $H_{\max / \min}^*(M, g, a)$; for the latter
case, however, we must calculate the weighted de Rham cohomology of each neighbourhood 
$\calU_\alpha$, which is the main goal of this section. We shall see that these 
are the same as the local intersection cohomology with respect to a certain
perversity function $\frakp$ depending on $a$, and this will prove the
equality of weighted de Rham and intersection cohomologies. 

By a slight abuse of notation, we denote the local cohomology of these sheaves by
$H^*_{{\max}/{\min}}(\calU_\alpha,g,a)$, respectively. If $\calU_\alpha$ is disjoint 
from all singular strata, then 
\[
H^j(\calU_\alpha,\calL_{{\max}/{\min},a}) = \left\{
\begin{array}{lcl}
\RR & \qquad & j = 0 \\
0 & \qquad & j > 0.
\end{array}
\right.
\]
Neither the weight function $x^a$ nor the metric $g$ play a role.
On the other hand, first note that 
\[
(\calL_{{\max}/{\min},a}(\calV \times C_1(F)), d_{{\max}/{\min},a})
= (L^2\Omega^*(\calV), d_{{\max}}) \widehat{\otimes} 
(L^2\Omega^*(C_1(F),g,a),d_{{\max}/{\min},a}).
\]
Furthermore, $\dim H^j(\calV) < \infty$ so the range of $d_{{\max}}$ on 
$L^2\Omega^*(\calV)$ is closed. Hence by Proposition~\ref{pr:kunneth} and
quasi-isometry invariance, 
\begin{eqnarray*}
H^j_{{\max}/{\min}}(\calV \times C_1(F),g,a) &=& 
\bigoplus_{\ell=0}^j 
H^\ell(\calV) \otimes H^{j-\ell}_{{\max}/{\min}}(C_1(F), dx^2 + x^2 \kappa ,a) \\
& = & H^j_{{\max}/{\min}}(C_1(F), dx^2 + x^2 \kappa ,a).
\end{eqnarray*}

We have now reduced the computation of weighted de Rham cohomology 
on $X$ to that of the truncated cone. For the next step we use
the following regularity result
\begin{lemma}
\[
\bigcap_{\ell=1}^\infty \calD(\Delta^\ell_{\Kr/\Fr},a) 
\subset \calA^*\Omega^*(C_1(F)).
\]
\label{le:conregdom}
\end{lemma}
\begin{proof}
The proof relies on the mapping properties of a parametrix for 
$\Delta_{{\Kr}/{\Fr},a}$; this parametrix is an element $G\in x^2\Psi_b^{-2,{\mathcal E}}
(M,\Omega^*)$, the calculus of $b$-pseudodifferential operators on $M$, such
that $G\Delta = I - R$ for some residual element $R$ which
satisfies $R:L^2\Omega^* \to \calA^*\Omega^*$. (Both $G$ and $R$
depend on $a$ and the choice of extension.) Suppose that 
$u \in \calD(\Delta_{\max,a})$, so in particular
$u, \Delta u \in x^a L^2\Omega^*$. Then $u = Gf + R u$, hence 
$u \in x^{a+2} H^2_b\Omega^* + \calA^*$. Induction on $\ell$ gives the result. 
We refer to \cite{Ma-edge} and \cite{Le} for more details.
\end{proof}

Now we prove three lemmas taken from \cite{Ch2}. In the following,
$g = dx^2 + x^2 \kappa$ on $C_1(F)$. 

\begin{proposition}
Let $r: C_1(F)\setminus \{0\} \to F$ be the canonical retraction map. Then
\[
r^*: L^2\Omega^k(F) \cap \ker d \longrightarrow \calD(d_{{\max},a}) \cap
x^a L^2\Omega^k(C_1(F))
\]
is well-defined and bounded if and only if $k < \frac{f+1}{2}-a$.
\label{prop:pullbacks}
\end{proposition}
\begin{proof}
Let $\alpha \in L^2\Omega^k(F)$; then
\[
\int_{C_1(F)}|\alpha|^2_g x^{-2a}\, dV_g = ||\alpha||^2_{L^2(F)} \int_0^1 x^{f-2k-2a} \, dx,
\]
and this is finite if and only if $k < \frac{f+1}{2}-a$. Since $d r^* \alpha = r^* d\alpha = 0$,
the image of $r^*$ lies in $\calD(d_{{\max},a})$.
\end{proof}

\begin{proposition}
There exists a $K>0$ such that for any $\omega = \alpha + dx \wedge \beta \in x^a L^2\Omega^k(C(F))$ 
there exists some $c \in (1/2, 1)$ for which
\[
||\alpha(c)||^2_{L^2(F)} \leq K ||\omega||^2_{x^aL^2(C_1(F))}.
\]
\label{prop:bound}
\end{proposition}
\begin{proof}
If not, then for any $N>0$, there exists some (nontrivial) $\omega\in x^aL^2\Omega^k(C_1(F))$ 
such that $N ||\omega||^2_{x^a L^2(C_1(F))} < ||\alpha(x)||^2_{L^2(F)}$ for all
$x \in (1/2,1)$. However, if this is the case, then for this $\omega$,
\[
||\omega||^2_{x^a L^2(C_1(F))} \geq 
||\alpha||^2_{x^a L^2(C_1(F))} \geq \int_{1/2}^1 ||\alpha(x)||^2_{L^2(F)} x^{f-2k-2a} \,dx 
\]
\[
> N ||\omega||^2_{x^a L^2(C_1(F))} \int_{1/2}^1 x^{f-2k-2a} \, dx. 
\]
This is a contradiction when $N$ is large.
\end{proof}

\begin{proposition} If $k < \frac{f+3}{2}-a$, then for any $c \in (1/2,1)$, the map
\[
x^a L^2 \Omega^k(C_1(F)) \ni\omega = \alpha + dx \wedge \beta \longrightarrow 
K_c(\omega) = \int_c^x \beta(s) \, ds \in x^a L^2 \Omega^{k-1}(C_1(F))
\]
is bounded.
\label{prop:K}
\end{proposition}
\begin{proof} First,
\begin{equation}
||K_c \omega ||^2_{x^a L^2(C_1(F))} = \int_0^1 \int_F \left|
\int_c^x \beta(s) \, ds\right|^2_\kappa \, x^{f-2k+2-2a} \, dx \, dV_\kappa.
\label{eq:Ka}
\end{equation}
Next, for any real number $b$, 
\[
\left|\int_c^x \beta(s) \, ds\right|^2_{\kappa} \leq 
\left( \int_c^x |\beta(s)|_\kappa \, ds \right)^2 
= \left( \int_c^x s^{-b} | s^b \beta(s)|_\kappa \, ds \right)^2
\]
\[
\leq \int_c^x s^{-2b} \, ds \, \int_c^x |\beta(s)|^2_\kappa s^{2b} \, ds
= \left\{ 
\begin{array}{ll}
\frac{x^{1-2b} - c^{1-2b}}{1-2b} \int_c^x |\beta(s)|^2_\kappa s^{2b} \, ds & b \neq 1/2 \\
(\log (x) - \log (c)) \int_c^x |\beta(s)|^2_\kappa s^{2b} \, ds & b = 1/2.
\end{array}
\right.
\]
Setting $2b = f-2k + 2-2a$ and using this in (\ref{eq:Ka}) gives
\[
||K_c \omega||^2_{x^a L^2(C_1(F))} 
\]
\[
\leq \left\{ 
\begin{array}{ll}
\int_0^1 \frac{x - c(x/c)^{f-2k + 2-2a}}{1-(f-2k+2-2a)} \int_c^x
|\beta(s)|^2_\kappa s^{f-2k+2-2a}\, ds \, dV_\kappa \, dx & k \neq  \frac{f+1}{2} -a\\
\int_0^1 x(\ln(x) - \ln(c)) \int_c^x |\beta(s)|^2_\kappa s^{f-2k+2-2a} \, ds \, dV_\kappa \, dr
&  k =  \frac{f+1}{2} -a
\end{array}
\right. 
\]
\[
\leq
\left\{ 
\begin{array}{ll}
\left(\int_0^1 \frac{x - c(x/c)^{f-2k + 2-2a}}{1-(f-2k+2-2a)} \, dx\, \right)
|| \beta||^2_{x^a L^2(C_1(F))}
& k \neq  \frac{f+1}{2} -a\\
\left(\int_0^1 x(\ln(x) - \ln(c)) \, dx\right) ||\beta||^2_{x^a L^2(C_1(F))}&  k =  \frac{f+1}{2} -a.
\end{array}
\right.
\]
Since $c$ is bounded away from $0$, both of these coefficients on the right
are uniformly bounded when $k \leq \frac{f+3}{2} -a$.
\end{proof}

\subsection{Poincar\'{e} lemma}
We now compute the weighted cohomologies of the truncated cone.
\begin{lemma} 
\[
H^k_{{\max}}(C_1(F),g,a) = \left\{
\begin{array}{lcl}
H^k(F) & \qquad & k < (f+1)/2 - a \\
0 & \qquad & k \geq  (f+1)/2 - a,
\end{array}
\right. 
\]
and 
\[
H_{{\min}}^k(C_1(F),g,a) = \left\{
\begin{array}{lcl}
H^k(F) & \qquad & k \leq (f-1)/2 -a  \\
0 & \qquad & k >  (f-1)/2 -a.
\end{array}
\right. 
\]
\label{lem:poincare}
\end{lemma}
\begin{proof} By Lemma \ref{le:conregdom}, we may work exclusively with conormal forms.
First let $k < \frac{f+1}{2} -a$. We wish to define a map
\[
R: H^k_{{\max}}(C_1(F),g,a) \longrightarrow H^k(F),
\]
and show that it is an isomorphism. Let $[\omega] \in H^k_{{\max}}(C_1(F),g,a)$, 
and choose a conormal representative $\omega = \alpha + dx \wedge \beta \in \maxdom$. 
For any $c \in (1/2,1)$, set $R([\omega]) = [\alpha(c)]$.  

To check that $R$ is defined independently of all choices, first note that
$d \omega = d_F \alpha + dx \wedge (\alpha' - d_F \beta) = 0$, so
$d_F \alpha(c) = 0$ for any $c$. Next, if $\tilde{\omega} = \tilde{\alpha} + 
dx \wedge \tilde{\beta}$ is another conormal representative of  $[\omega]$, 
then there exists $\eta = \mu + dx \wedge \nu \in x^a L^2\Omega^{k-1}(C_1(F))
\cap \calA^*$ with $\tilde{\omega} = \omega + d \eta$. This implies in particular
that $\tilde{\alpha}(c) = \alpha(c)+ d_F \mu(c)$, so $[\alpha(c)] = [\tilde{\alpha}(c)]$.
Similarly, $\int_c^{c'} \beta(s) \, ds \in L^2\Omega^{k-1}(F)$ and hence
$\alpha' = d_F \beta$ implies that
\[
d\int_c^{c'} \beta(s) \, ds  = \int_c^{c'} d_F \beta(s) \, ds = \int_c^{c'} \alpha'(s) \, ds
= \alpha(c') - \alpha(c);
\]
thus $[\alpha(c)] = [\alpha(c')]$. 

$R$ is certainly linear; it is bounded by Proposition \ref{prop:bound} and surjective 
by Proposition \ref{prop:pullbacks}, so we must only show that it is injective.  
Suppose $\omega = \alpha + dx \wedge \beta \in \calA^*$ and $R([\alpha])=[0]$. 
Then $\alpha(c) = d_F \eta$ for some $\eta \in \calC^\infty\Omega^{k-1}(F)$. By 
Proposition \ref{prop:pullbacks}, $r^*\eta \in \maxdom$, while Proposition \ref{prop:K} 
gives that $\int_c^x \beta(s) \, ds \in x^a L^2\Omega^{k-1}(C_1(F)) \cap \calA^*$. Thus
\[
d\left(\eta + \int_c^x \beta(s) \, ds\right) = 
d_F \eta + dx \wedge \beta(x) + \int_c^x d_F\beta(s)\, ds 
\]
\[
= \alpha(c) + dx \wedge \beta(x) + \int_c^x \alpha'(s)\, ds = \omega.
\]
Since $\omega \in  x^a L^2\Omega^k(C_1(F))$, this implies that
$\eta + \int_c^x \beta(s) \, ds \in \maxdom$, hence $\omega$ is exact in the 
maximal complex, i.e.\ $[\omega] = [0]$, as desired.

We next show that when $k \geq \frac{f+1}{2} -a$, any $[\omega] \in 
H^k_{{\max}}(C_1(F),g,a)$ is trivial. Thus, for any representative
$\omega = \alpha + dx \wedge \beta \in \maxdom \cap \calA^*$, we must
find a $(k-1)$-form $\eta \in \maxdom\cap \calA^*$ with $\omega = 
d_{{\max},a}\eta$. Assume $|\omega|_\kappa = O(x^p)$ for some $p$. 
The condition $\int_0^1 |\omega|^2_\kappa x^{f- 2k - 2a}\, dx dV_\kappa < \infty$ 
gives that $p > a+k-\frac{f+3}{2}$. Furthermore, by assumption, $a+k-\frac{f+3}{2} > -1$, 
so $K_0 (\omega) = \int_0^x \beta(s) \, ds$ is defined. Using $|\beta(s)|_\kappa = 
O(s^{p+1})$, we deduce that the integral
\[
||K_0 \omega||^2_{x^a L^2(C_1(F))} = \int_0^1 \int_F \left| \int_0^x \beta(s)\, ds
\right|_\kappa^2 x^{f-2k + 2-2a} \, dx dV_\kappa
\]
is finite.  Now, $\alpha \in \calA^* \cap x^a L^2 \Omega^{k}(C_1(F))$, so 
$|\alpha|_\kappa = O(x^q)$ where $2q + f - 2k - 2a > -1$, i.e.\ $q > 0$, so
\[
d(K_0(\omega)) = dx \wedge \beta(x) + \int_0^x \alpha'(s) \,ds = \omega(x),
\]
as desired.  We have now shown that $K_0(\omega) \in \maxdom$ and hence 
$[\omega] = 0$ in $H^k_{{\max}}(C_1(F),g,a)$. This completes the computation 
of $H^*_{{\max}}(C_1(F),g,a)$.

\medskip

The computation of $H^*_{{\min}}(C_1(F),g,a)$ proceeds identically when
$k \leq (f-1)/2 -a$ or $k \geq (f+3)/2-a$, but the remaining cases are
treated slightly differently.

\begin{claim}
If $k \in ((f-1)/2 - a, (f+1)/2) -a)$ and $\omega = \alpha + dx\wedge\beta
\in x^a L^2 \Omega^k(C_1(F)) \cap \calD(d_{{\min},a})\cap \calA^*$ is closed, 
then $\omega=d\eta$ for some $\eta \in \calD(d_{{\min},a})\cap \calA^*$. As a
consequence, $H_{\min}^k(C_1(F),g,a)= 0$.
\end{claim}

\begin{proof}
As above, $\alpha' = d_F\beta$, so if $c,c' \in (0,1)$, 
\[
\Pi_0(\alpha(c') - \alpha(c)) = \Pi_0 \int_c^{c'} \alpha'(s) \,ds = \Pi_0 d_{F}\int_c^{c'} 
\beta(x) \,dx = 0.
\]
Since $\omega \in \calD(d_{{\min},a})$, $\Pi_0(\alpha(x))=o(1)$, so $\Pi_0(\alpha(x))=0$ for 
all $x$. This gives that
\[
d\int_c^x \beta(s)\,ds = dx \wedge\beta + \alpha(x)-\alpha(c) 
= \omega - \alpha(c) = \omega - \Pi_{\perp}\alpha(c) = \alpha - d(r^* \eta)
\]
for some $\eta \in L^2 \Omega^{k-1}(F) \cap \calD(d)$. 
Thus $\alpha = d(\eta + \int_c^x \omega(s) \, ds)$.

As for its domain, it suffices by Lemma \ref{le:max=min} to show that 
$\eta + \int_c^x \beta(s) \, ds \in \calD(d_{{\max},a})$ since $k-1 < (f-1)/2 - a$.
But on the one hand, $d(\eta + \int_c^x \beta(s) \, ds)= \omega \in x^a L^2\Omega^{k}(C_1(F))$;
furthermore, $(\eta + \int_c^x \omega(s) \, ds) \in  L^2\Omega^{k-1}(C(F), g,a)$ 
since $r^*\eta \in x^a L^2\Omega^{k-1}(C_1(F))$, again because $k-1 < (f+1)/2-a$, 
and $\beta = o(x^{(2k-2+2a-f-1)/2})= o(x^{p})$ for some $p > -1$, so
$\int_c^x \beta(s)\, ds \in x^a L^2\Omega^{k-1}(C_1(F))$.
This proves the claim.
\end{proof}

\begin{claim}
The map $H_{\min}^k(C_1(F),g,a) \rightarrow H^k(F)$ is injective when
$k \in ((f+1)/2 - a, (f+3)/2 -a)$; hence, for $k$ in this range, 
$H_{\min}^k(C_1(F),g,a)\cong H_{max}^i(C(F),g,a) \cong 0$.
\end{claim}

\begin{proof}
If $\omega = \alpha + dx \wedge \beta \in \calA^*$ represents a class in 
$H_{\min}^k(C_1(F),g,a)$, then we already know that $\omega = d\eta$ for 
some $\eta \in\calD(d_{{\max},a}) \cap \calA^*$. So we must show that 
we can arrange for $\eta$ to lie in $\calD(d_{{\min},a})$ as well. 
Since $k-1 \in ((f-1)/2 + a, (f+1)/2 -a)$, $\beta(x) = o(x^{p})$ for some $\nu > -1$,
and so
\[
\Pi_0 \int_c^0 \beta(s)ds
\]
is defined. We can choose $\zeta$ which solves $\alpha(c) = d_{F}\zeta$ 
by specifying that $\Pi_0(\zeta) = -\Pi_0  \int_c^0 \beta(s)ds$. Now let 
$\eta = \zeta + \int_0^x \beta(s)ds$.  As before, $\eta \in\calA^*$ and 
$d\eta = \omega$.  Furthermore, $\Pi_0(\eta(0))=0$ and $\del_x(\Pi_0(\eta))= 
\Pi_0(\beta(x))=o(x^{(2(k-1)+2a-f-1)/2})=o(x^{p})$ for $p>-1$, so 
$\Pi_0(\eta) = o(1)$. This means that $\eta \in \calD(d_{{\min},a})$.
\end{proof}

\end{proof}

This completes the calculation of the local cohomology for the sheaves 
$\calL_{{\max}/{\min},a}$ on $X$. By Proposition~\ref{pr:shchar}, we now
obtain one of our main results:
\begin{theorem}
If $(M,g)$ is a manifold with an incomplete edge metric and $X$ is the associated
stratified space, then 
\[
H^*_{\max}(M,g,a) =
\left\{ 
\begin{array}{lll}
I\!H^*_{\overline{\frakm} + \ll a-1 \gg} (X,B) & \qquad & f \mbox{ odd} \\
I\!H^*_{\overline{\frakm}  + \ll a-1/2 \gg } (X,B) & \qquad & f \mbox{ even} 
\end{array}
\right.
\]
and
\[
H^*_{\min}(M,g,a) = 
\left\{ 
\begin{array}{lll}
I\!H^*_{\underline{\frakm} + <a>} (X,B) & \qquad & f \mbox{ odd} \\
I\!H^*_{\underline{\frakm} + <a-1/2>} (X,B) & \qquad & f \mbox{ even} 
\end{array}
\right.;
\]
\noindent
here $\ll t \gg $ denotes the least integer strictly greater than $t$ and $<t>$ denotes the 
least integer greater than or equal to $t$.
\label{th:weightedco}
\end{theorem}

There are two important special cases which we single out:
\begin{corollary}
The maximal and minimal de Rham cohomologies at weight zero correspond
to upper and lower middle perversity intersection cohomology.
\begin{equation}
\begin{array}{rcl}
H^*_{\max}(M,g,0) &= &  I\!H^*_{\overline{\frakm}}(X) \\
H^*_{\min}(M,g,0) & = &  I\!H^*_{\underline{\frakm}}(X)
\end{array}
\label{eq:wt0}
\end{equation}
Moreover, when $f$ is even, the maximal and minimal de Rham
cohomologies at weights $\pm 1/2$ coincide, and again correspond
to upper and lower middle perversity intersection cohomology.
\begin{equation}
\begin{array}{rcl}
H_{{\max}/{\min}}^*(M,g,-1/2) &=& I\!H^*_{\overline{\frakm}}(X)\\
H_{{\max}/{\min}}^*(M,g,1/2)  & = & I\!H^*_{\underline{\frakm}}(X).
\end{array}
\end{equation}
\label{cor:hmaxmin}
\end{corollary}

{}From \S 2, these weighted de Rham cohomology spaces are identified
with the nullspaces of the associated (absolute and relative)
Laplacians, and we conclude that the nullspaces of $\Delta_{\Kr/\Fr,a}$ 
are (finite dimensional and) identified with particular intersection
cohomology spaces. Note that we are {\it not} asserting anything about the
nullspaces of the  `ordinary' Laplacians $\Delta_{\Kr/\Fr,0}$ on the weighted 
spaces $x^a L^2\Omega^*$;  indeed, these nullspaces are either infinite dimensional
when $a \ll 0$ or vanish identically when $a \gg 0$, cf. \cite{Ma-edge}.

\section{Elliptic edge operators and minimal Hodge cohomology} 
To proceed further in the study of these weighted de Rham complexes, we must use elliptic 
methods. More specifically, we still wish to study the question of when there is a unique 
closed extension for $d$ on $x^a L^2 \Omega^*(M)$ for incomplete edge metrics, and we also 
wish to compute the minimal Hodge cohomology. We shall study both of these questions 
using the formally symmetric operator $D_a = d+\delta_a$ on $x^a L^2 \Omega^*(M)$. 
The proper context for this analysis is the calculus of pseudodifferential edge
operators, and in the next subsection we review the generalities of this theory. After 
that we show how it applies to the specific problems at hand.

\subsection{Edge operators}
We now review the general theory of elliptic edge operators. This is the correct 
context to study $D_a$ for an incomplete edge metric (and also the corresponding 
operator for a complete edge metric). This theory is developed fully in \cite{Ma-edge},
and we refer there for more details

Fix a local coordinate system $(y_1, \ldots, y_b)$ on $B$ and $(z_1, \ldots, z_f)$
on $F$, so that $w = (x,y,z)$ is a local coordinate system in some
neighbourhood of a singular stratum in $X$. By including the hypersurface
$\{x=0\}$, we are blowing up $B$ in $X$; the resulting manifold
with boundary is denoted $\overline{M}$ and its interior is denoted $M$. 

A differential operator $L$ on $M$ is called an edge operator of order $m$ 
if it can be expressed in the form 
\begin{equation}
L = \sum_{j+|\alpha| + |\beta| \leq m} a_{j,\alpha}(x,y,z)(x \, \del_{x})^j 
(x \, \del_{y})^\alpha \del_z^\beta,
\label{eq:unifdeg}
\end{equation}
where the (scalar or matrix-valued) coefficients are bounded. We shall assume that 
these coefficients are smooth in these variables, down to $x=0$. 
For example, if $g$ is a {\it complete} edge metric, then the scalar 
or Hodge Laplacian is an operator of this type; similarly, if $g$ is
an {\it incomplete} edge metric, then its Laplacian is of the form $x^{-2}L$, 
where $L$ is an edge operator of order $2$.

\subsubsection{Ellipticity and model operators}
There is a well-defined symbol in this setting:
\[
\sigma(L)(x,y,z;\xi,\eta,\zeta) : = \sum_{j+|\alpha| + |\beta| = m} 
a_{j,\alpha}(x , y, z)\, \xi^j \eta^\alpha  \zeta^\beta,
\]
and we say that $L$ is elliptic in the edge calculus provided $\sigma(L)(x,y,z;\xi,\eta,\zeta)$ is 
invertible when $(\xi,\eta,\zeta) \neq 0$. 

Ellipticity alone does not guarantee that $L$ is Fredholm between appropriate function spaces; 
one must also require that certain model operators for $L$ also be invertible. 
There are two such operators: 
\begin{itemize}
\item The {\it normal operator} of $L$ is defined by 
\[
N(L) : = \sum_{j+|\alpha| + |\beta| \leq m} a_{j,\alpha,\beta}(0,y,z)(s \del_{s})^j 
(s \del_{u})^\alpha \del_z^\beta\, \qquad (s,u) \in \RR^+ \times \RR^b ;
\]
here $y\in B$ enters only parametrically and the operator acts on  functions on 
$\RR^+ \times \RR^b \times F$. This operator can be regarded as $L$ with its 
coefficients frozen (in an appropriate sense) at $x=0$, acting on functions (or 
sections of an appropriate bundle) on the space $\RR^+_s \times \RR^b_u \times F_z$. 
\item The {\it indicial operator} of $L$ is defined by
\[
I(L) : = \sum_{j +|\beta| \leq m} a_{j,0,\beta}(0,y,z) (s\del_s)^j \del_z^\beta.
\]
\end{itemize}

For example, the normal and indicial operators associated to the scalar Laplacian
for the complete edge metric $x^{-2}(dx^2 + h) + \kappa$ are 
\[
N(\Delta_{g}) = s^2 \, \del_s^2 + (1-b) \, s \, \del_s + s^2 \Delta_u + \Delta_{\kappa}, \qquad
I(\Delta_{g}) = s^2 \, \del_s^2 + (1-b) \, s \, \del_s + \Delta_{\kappa}.
\]

The indicial operator captures some fundamental invariants associated to $L$: 
\begin{definition} 
The number $\gamma \in \CC$ is said to be an indicial root of $L$ at $y_0 \in B$
if there exists a function $v(z)$ on $F$ such that (in local coordinates where 
$y_0$ corresponds to $y=0$)
\[
I(L)_{y_0}(s^\gamma v(z)) = \left(\sum_{j +|\beta| \leq m} a_{j,0,\beta}(0,0,z)
(s\del_s)^j\right) s^\gamma v(z) = \calO(s^{\gamma+1}). 
\]
\end{definition}
Indicial roots may often be calculated in terms of eigenvalues for an induced elliptic 
operator on the fibre $F_{y_0}$, and might depend on $y_0 \in B$. 

The operator $L$ acts naturally on weighted Sobolev spaces. Let $M$ be a manifold 
with {\it complete} edge metric $G$.  For $\ell \in \NN$ and $\delta \in \RR$, define 
\[
x^\delta H^\ell_{e}(M) = \{u = x^\delta v: (x\del_x)^j(x\del_y)^\alpha\del_z^\beta v
\in L^2(M,dV_G)\ \forall\ j+|\alpha|+|\beta| \leq \ell \}.
\]
(By interpolation and duality, these spaces can be defined for any $\ell \in \RR$.)
Clearly, if $L$ is any edge operator of order $m$, then 
\begin{equation}
L: x^{\delta}H^{\ell+m}_e(M) \longrightarrow x^\delta H^\ell_e(M)
\label{eq:mapwss}
\end{equation}
for any $\delta, \ell$. Further hypotheses, beyond the ellipticity of $L$ are required
to ensure that this mapping is well-behaved. 

The first instance of this is that the indicial roots of $L$ yield weights $\delta$ 
for which (\ref{eq:mapwss}) does not have closed range; these are precisely the 
weights $\delta$ for which an indicial root $\gamma$ `just fails' to lie in $x^\delta L^2$ 
near $x=0$, i.e.\ where $x^\gamma \in x^{\delta - \epsilon}L^2$ for any $\epsilon > 0$ 
but $x^\gamma \notin x^\delta L^2$. We denote this critical weight $\delta$ associated
to a given indicial root $\gamma$ as $\delta(\gamma)$. With respect to the measure $dxdydz$, 
$\delta(\gamma) := \mbox{Re}\,\gamma + 1/2$. However, the measure appearing in our application 
below is $x^{f-2a}\, dxdydz$, and $x^\gamma \in x^\delta L^2(x^{f-2a}\, dxdydz)$ near $x=0$ 
if and only if $\delta < \gamma + (f+1)/2 - a$, so that we shall define $\delta(\gamma) 
= \mbox{Re}\, \gamma + (f+1)/2 - a$. 

Even when $\delta$ is not equal to one of these critical values, the behaviour of
the normal operator at weight $\delta$ plays another very important role. 

\begin{proposition}
Let $L$ be an elliptic differential edge operator of order $m$. Fix $\delta$ such that
$\delta \neq \delta(\gamma(b))$ for any indicial root $\gamma(b)$, $b \in B$. Suppose 
also that $N(L): s^\delta H^m_e \to s^\delta L^2$ is surjective (for all $b \in B$). 
Then (\ref{eq:mapwss}) is essentially surjective, in the sense that its range
is closed and of finite codimension.  On the other hand, if $N(L)$ is
injective on $s^\delta L^2$, then any element of the nullspace of $L$
is necessarily conormal. 
\end{proposition}

There are many more refined statements one can make about the mapping
properties of $L$. For later applications, we state only one very special
result. We shall restrict to a special setting, which is what arises in our 
applications below. The hypothesis that the the normal operator $N(L)$ is surjective 
is equivalent to the injectivity of the normal operator for the adjoint $L^*$. 
This adjoint depends on the choice of measure, and we shall assume (as
in our applications) that the adjoint of $L$ on $x^\delta L^2$ corresponds 
to the same operator $L$ on a `dually weighted' space $x^{\delta^*}L^2$,
for some $\delta^* > \delta$. The fact that $L$ has closed range implies the 
existence of a generalized inverse $G: x^\delta L^2 \to x^\delta H^m_e$ which 
satisfies $LG = I - P$, where $P$ is the orthogonal projector onto the cokernel.
By duality, elements of this cokernel are identified with elements of
the nullspace of $L$ on $x^{\delta^*}L^2$, and by the result above, 
these are conormal. 

\begin{proposition} 
Let $L$ satisfy the special assumptions of the preceding paragraph.
Suppose furthermore that the interval $(\delta,\delta^*)$ contains 
a finite set of indicial roots $\gamma_j$, $j=1,\ldots, N$, 
all of which are constant in $b \in B$.  Let $f \in x^{\delta^*}L^2
\cap \calA^*$. Then $u = Gf$ satisfies $Lu = f - \phi$ where
$\phi \in x^{\delta^*}L^2\cap \calA^*$, $L\phi = 0$, and
$u = \sum_{j=1}^N u_j(y,z) x^{\gamma_j} + v$; where each $u_j(y,z)\in \calC^\infty$ 
solves the indicial equation $I(L)(s^{\gamma_j} u_j(y,z)) = 0$, and the error term 
$v \in \calA^* \cap x^{\delta^*}L^2$.
\end{proposition}

The proofs rely on the construction of a pseudodifferential parametrix $G$ for $L$, 
depending on $\delta$. This is an element in the calculus of pseudodifferential 
edge operators $\Psi^*_e(M)$. We do not define this calculus here, but remark 
only that these operators are described by specifying the precise asymptotic 
behaviour of their Schwartz kernels, near the diagonal and also near the boundaries 
and corner of $M \times M$. We refer as before to \cite{Ma-edge}. (We should note 
also that the results stated here are slightly more general than what is written 
explicitly in that source because we are allowing the possibility of variable
indicial roots outside the critical interval; however, these can be derived 
easily from the same techniques.)

\subsection{Edge analysis of $D_a$}
We now proceed to apply the methods of the last subsection to the analysis
of $D_a$. The first tasks are the calculation of the indicial roots of 
$D_a$ and the analysis of the normal operator $N(xD_a)$.

\subsubsection{Indicial roots of $D_a$}

For simplicity, first consider the calculation of the indicial roots for $D_a$
for the metric $g = dx^2 + x^2 \kappa$ on the cone $C_1(F)$. To do this, we decompose
this operator as much as possible. Thus, first regard $D_a$ as a $2\times 2$ matrix 
acting on pairs $(\alpha,\beta) \leftrightarrow \alpha + dx \wedge \beta$, where
$\alpha(x), \beta(x) \in \Omega^*(F)$ for each $x$; normalize by writing the $k$-form 
part of $\alpha$ as $x^k \alpha_k$, and similarly for $\beta$. A short calculation then
shows that, acting on pairs $(\alpha_k,\beta_k)$, 
\[
\left. I(xD_a) \right|_{\Omega^k \oplus \Omega^k} = \left(\begin{array}{ll}
\frac{1}{x} D_F & -\del_x - \frac{f-k-2a}{x} \\
\del_x + \frac{k}{x} & -\frac{1}{x}D_F \end{array}\right);
\]
the full indicial operator is the direct sum over $k$ of these matrices. 
Similarly, the indicial family is the direct sum of matrices 
\begin{equation}
I(xD_a)_k(\gamma) = \left(\begin{array}{ll}D_F & -\gamma - (f-k+2a) \\
\gamma + k & -D_F\end{array}\right)
\label{eq:indfamda}
\end{equation}
This can be reduced further using the eigendecomposition for $\Delta_F$. In particular,
we see that the operator in (\ref{eq:indfamda}) is noninvertible if and only if for some
eigenvalue $\lambda^2$ for $\Delta_F$,  
\[
\left(\begin{array}{ll}
\lambda & -\gamma - (f-k-2a) \\
\gamma + k & -\lambda \end{array}\right)
\]
is singular, or equivalently
\[
\gamma^2 + (f - 2a)\gamma + k(f-k-2a) - \lambda^2 = 0.
\]
Hence the indicial roots come in pairs:
\begin{equation}
\gamma_{\lambda,k}^{\pm} = a -\frac{f}{2} \pm \frac12 \left[(f-2a - 2k)^2 + 4\lambda^2\right]^{1/2}.
\label{eq:enumir}
\end{equation}

The extension of these calculations to incomplete edge metrics requires only some mild 
alterations. Write $g = dx^2 + \tilde{g}(x)$, where $\tilde{g}(x) = h + x^2 \kappa$ is 
a degenerating family of metrics on $Y = \del M$. When $\alpha$ is a $(p,q)$-form on $Y$, 
its pointwise norm satisfies
\[
|\alpha|_{\tilde{g}(x)} = x^{-q} |\alpha|_{\tilde{g}(1)}.
\]
Denote by $D_Y^x$ the operator $D$ on $Y$ for the metric $G_x$. By
Proposition 3, with respect to the metric $\tilde{g}(1)$ on $Y$,
\[
d_Y = d_F + \tilde{d}_B - \II + {\rm R}, \qquad 
\delta_Y = \delta_F + (\tilde{d}_B)^* - \II^* + {\rm R}^*;
\]
hence, $D_Y = D_F + \overline{D}_B + \overline{R}$ where
\[
\overline{D}_B = \tilde{d}_B + (\tilde{d}_B)^* - \II - \II^*, \qquad
\mbox{and}\quad \overline{R}={\rm R} + {\rm R}^*.
\]
A quick review of the definitions shows that 
\[
D_Y^x = \frac{1}{x}D_F + \overline{D}_B + x \overline{R},
\]
where all the components on the right are the corresponding operators at $x=1$. 

Now let $Z \to B$ be the bundle with fibre $C_1(F)$ obtained from $Y \to B$, with 
metric $g = dx^2 + \tilde{g}(x)$. Decompose any form $\omega = \alpha + dx \wedge \beta$ 
on $Z$ as $\alpha = \sum x^k \alpha_k$, where $\alpha_k$ is of type $(*,k)$ on $Y$, 
and similarly for $\beta$. Thus the pointwise norms satisfy
$|\omega|^2_g = \sum \left(|\alpha_k|_{\tilde{g}(1)}^2 + |\beta_k|_{\tilde{g}(1)}^2\right)$. 
In terms of these decompositions and normalizations, the restriction of the
operator $D_a$ to pairs of $(*,k)$-forms on $Z$ is given by
\begin{equation}
D_a = \left(\begin{array}{ll}
\frac{1}{x} D_F + \overline{D}_B + x \overline{R} & -\del_x - \frac{f-k-2a}{x} \\
\del_x + \frac{k}{x} & -\left(\frac{1}{x}D_F + \overline{D}_B + x\overline{R}\right)
\end{array}\right).
\label{eq:Daz}
\end{equation}

{}From this expression, we see that neither $\overline{D}_B$ nor $\overline{R}$ 
appear in the indicial operator
$I(xD_a)$. Hence the computation of the indicial roots is exactly the same as in the 
conic case; in other words, all indicial roots are of the form (\ref{eq:enumir}). Note,
however, that the eigenvalues $\lambda^2$ may depend on $b \in B$, hence the same may
be true of these indicial roots. 

Notice that if $\omega \in \Omega^k$ and $|\omega|_g \sim x^\gamma$, then
$\omega \in x^a L^2\Omega^k$ (near $x=0$) if and only if $\gamma > a - f/2$.
The indicial roots which lie near to (and above) this `$x^a L^2$ cutoff' 
are the ones which cause the difference between minimal and maximal domains.
We explain this later, but for now record the 
\begin{corollary} The indicial roots of the operator $D_a$ contained in the interval
$(a-(f+1)/2, a - (f-1)/2)$ correspond to the eigenvalues $\lambda^2$ of $\Delta_F$ 
on $k$-forms such that $(f-2a + 2k)^2 + 4\lambda^2 < 1$. In order for this condition
to be nonvacuous, it is necessary that $k \in ((f-1)/2 -a ,(f+1)/2 - a)$. Note in particular
that $\gamma_{0,k}^{\pm} = a-f/2 \pm (f/2 - a - k) = -k, k+2a-f \in (a-(f+1)/2, a - (f-1)/2)$
precisely when $k \in ((f-1)/2 - a,(f+1)/2-a)$ and $H^k(F) \neq 0$. 
\end{corollary}

\subsubsection{The normal operator of $D_a$}

According to the discussion in the final paragraphs of \S 4.1, we must also
study the mapping properties of the normal operator $N(xD_a)$. Before doing
so, we address some `duality' issues. The main point is that $D_a$ is formally 
symmetric on $x^a L^2\Omega^*$, hence the adjoint of the (closed range) operator
\begin{equation}
D_a: x^a L^2 \Omega^* \longrightarrow x^{a-1}L^2\Omega^*
\label{eq:Daa-1}
\end{equation}
is identified with
\begin{equation}
D_a: x^{a+1} L^2 \Omega^* \longrightarrow x^{a}L^2\Omega^*.
\label{eq:Da+1a}
\end{equation}
Recall also that, according to the computations of the preceding subsection, there are at 
most two indicial roots $\gamma_0^\pm$ in the interval $(a-(f+1)/2, a-(f-1)/2)$, and these 
are symmetric around the midpoint $a-f/2$. We shall assume that the metric $g$ is
such that $\Delta_F$ has no small nonzero eigenvalues, so that no other indicial roots
intersect the closed interval $[a-(f+1)/2, a-(f-1)/2]$.

By (\ref{eq:Daz}), 
\begin{equation}
N(xD_a) = s D_{C(F),a} + s D_{\RR^b},
\label{eq:normop}
\end{equation}
where the first operator on the right is the analogous weighted operator on the
complete cone $C(F)$ and the second is on Euclidean space. In fact, the identification
of $s^{-1}N(xD_a)$ with $D_{C(F),a} + D_{\RR^b}$ may also be seen by naturality, since
the operator on the left must equal, at $b \in B$, the Hodge-de Rham operator on
$\RR^+_s \times \RR^b_u \times F$ with respect to $g_b = ds^2 + s^2 |du|^2 + \kappa_b$.

\begin{proposition} 
\[
N(xD_a): s^{a+1} L^2\Omega^* \longrightarrow s^{a+1} L^2 \Omega^*
\]
is injective, and hence (\ref{eq:Da+1a}) has a finite dimensional
nullspace consisting of conormal forms
\label{pr:injnDa}
\end{proposition}

\begin{proof} Suppose that $\omega \in s^\alpha L^2\Omega^*$ is in the
nullspace of this operator. Take the Fourier transform in the $u$ direction;
denoting the dual variable by $\eta$, then whenever $\eta \neq 0$ we can rescale, 
setting $t = s|\eta|$, $\hat{\eta} = \eta/|\eta|$. Then
\[
\left(D_{C(F),a} + i \, \mbox{cl}\,(\hat{\eta})\right) \hat{\omega} = 0,
\]
where $\mbox{cl}\,(\hat{\eta})$ is Clifford multiplication 
$\hat{\eta}\wedge \cdot + \iota(\hat{\eta}) \cdot$. Apply $D_{C(F),a} + 
i \, \mbox{cl}\,(\hat{\eta})$ to this equation to deduce that
\[
\left(\Delta_{C(F),a} + 1\right) \hat{\omega} = 0.
\]
It is not hard to show, cf.\ \cite{Ma-edge}, that any solution of this
equation either grows or decays exponentially as $t \to \infty$, and the
$L^2$ hypothesis prohibits the former. Furthermore, solutions are polyhomogeneous
as $t \to 0$, and hence decay at some indicial weight $t^\gamma$ 
with $\gamma > a - (f-1)/2$. Hence both $N(d)\hat{\omega}$ and $N(\delta_a) \hat{\omega}$ 
decay like $t^{\gamma - 1}$, and in particular are still in $t^a L^2$. 
This means we can integrate by parts to obtain
\[
0 = \langle (\Delta_{C(F),a} + 1)\hat{\omega},\hat{\omega}\rangle_a = 
||N(d)\hat{\omega}||_a^2 + ||N(\delta_a)\hat{\omega}||_a^2 + ||\hat{\omega}||_a^2;
\]
all boundary terms vanish. This gives $\hat{\omega} = 0$, as desired.

When $\eta = 0$, the problem reduces to showing that $\Delta_{C(F),a}$
has no nullspace in $t^{a+1}L^2$ on the entire cone $C(F)$, which is
even more easily verified to be true (e.g.\ by separation of variables).
\end{proof}

Following the discussion from the end of \S 4.1, if $(a-(f+1)/2, a-(f-1)/2)$
contains no indicial roots for $D_a$, $N(xD_a)$ is injective on $s^a L^2\Omega^*$.

There are no forms in the nullspace of the normal operator which lie in $t^a L^2\Omega^* 
\cap \calD(d_{{\max},a})\cap \calD(\delta_{{\min},a})$ or $t^a L^2\Omega^* \cap \calD
(d_{{\min},a})\cap \calD(\delta_{{\max},a})$. From this it is possible to show that 
$\omega \in \calH^*_{\Kr/\Fr}(M,g,a) \subset \calA^*\Omega^*$, i.e.\ such harmonic 
forms are conormal. Unfortunately, these results rely on a slightly more elaborate parametrix 
construction than is available in the literature, so at present we are only asserting
this informally. At the end of the next subsection, 
however, we show that forms in $\calH^*_{{\min}}(M,g,a)$ are conormal. 

\subsection{Closed extensions of $D_a$}
It is a general fact that closed extensions of $d$ on $x^a L^2\Omega^*$ are in
bijective correspondence with the self-adjoint extensions of $D_a$ on this space. 
To see this, first note that if $\overline{d}$ is any closed extension of $d$, and if 
$\overline{d_a}^*$ is its Hilbert space adjoint, then $\overline{D_a}= \overline{d} + 
\overline{d_a}^*$ is a self-adjoint extension of $D_a$. Conversely, any self-adjoint 
extension of $D_a$ determines an associated closed extension for $d$, 
cf.\ \cite[Lemma 2.3]{BL}. We summarize this in the
\begin{proposition}
If $d$ has more than one closed extension on $x^a L^2 \Omega^*(M,g)$, then $D_a$ has 
more than one self-adjoint extension on this space; equivalently, if $D_a$
is essentially self-adjoint, then $d_{{\max},a} = d_{{\min},a}$.
\label{pr:uce}
\end{proposition}
Notice that since $d$ always has closed extensions, $D_a$ always has self-adjoint 
extensions. However, if $D_a$ is not essentially self-adjoint, then it will have
closed extensions which are not necessarily self-adjoint; the relationship between 
these and the closed extensions of $d$ is somewhat more complicated, and we shall
not attempt to describe it.

\begin{theorem}
The symmetric operator $D_a$ is essentially self-adjoint on $x^a L^2\Omega^*(M,g)$ 
if and only if there exists no indicial root for $D_a$ in the interval $(a - (f+1)/2, 
a-(f-1)/2)$. As explained earlier, this is equivalent to the nonexistence of small 
eigenvalues $\lambda^2$ for $\Delta_F$ on $k$-forms such that $k \in (f/2 - a - 
\frac12\sqrt{1-4\lambda^2}), f/2 - a + \frac12\sqrt{1-4\lambda^2})$. 
In this case $d$ also has a unique closed extension. 
\end{theorem}

\begin{proof}
Assume that there are no small eigenvalues, as described in the statement of the theorem. 
Fix a parametrix $G$ for $D_a$ relative to the space $x^a L^2\Omega^*$. This is an 
element of order $-1$ in the edge calculus. If $\omega \in \calD(D_{{\max},a})$, 
then $f = D_a \omega \in x^a L^2\Omega^*(M,g)$, and applying $G$ gives that in fact 
$\omega \in x^{a+1}H^1_e\Omega^*$. We recall that in general, $Gf$ would
be the sum of two terms, the first corresponding to these small indicial roots and the
second an error term in $x^{a+1}H^1_e\Omega^*$, but by our hypothesis, the former 
of these is absent. It is now straightforward to check that $\omega$ may be smoothly 
approximated in the $D_a$-graph norm, i.e.\ that there exists a sequence of smooth 
compactly supported forms $\phi_j$ such that $\phi_j \to \omega$, $D_a \phi_j \to 
D_a \omega$ in $x^a L^2\Omega^*$. This shows that $\omega \in \calD(D_{{\min},a})$. 

Conversely, if there do exist indicial roots in the critical range, then these may
be used to construct nontrivial elements in $\calD(D_{{\max},a})\setminus
\calD(D_{{\min},a})$, and by the general abstract theory, there will be more than 
one self-adjoint extension of $D_a$. Because this is not central to our discussion,
we leave details to the reader (and refer to \cite{GilMen} for a thorough
discussion of the conic case). 
\end{proof}

In the conic case (when $a=0$) this result is due to Cheeger \cite{Ch1}; cf.\ 
also \cite{Le}. The analysis needed in that case is simpler than the edge 
analysis used here, though this is not apparent `on the surface'. 

We turn now to a description of $\calD(D_{{\max},a})$ in the more general
case where this vanishing condition is no longer satisfied. By definition, 
this domain is the set of all $\omega \in x^a L^2\Omega^*$ such that 
$D_a \omega \in x^a L^2\Omega^*$. We have already remarked that the mapping 
(\ref{eq:Daa-1}) has closed range, which we denote ${\mathcal R}_a$. 
Its cokernel is finite dimensional since, 
by Proposition (\ref{pr:injnDa}), the adjoint mapping (\ref{eq:Da+1a}) has 
a finite dimensional nullspace. Hence we can choose a generalized inverse
\[
G_a: x^{a-1}L^2\Omega^*(M,g) \longrightarrow x^a L^2\Omega^*(M,g);
\]
this is a pseudodifferential edge operator of order $-1$ which satisfies 
$D_a G_a = I - P$ on $x^{a-1}L^2\Omega^*$, where $P$ is the orthogonal projector 
onto the cokernel. The condition  $f \in {\mathcal R}_a$  is equivalent to 
$\langle f ,\gamma  \rangle = 0$ for all $\gamma$ in the nullspace of 
(\ref{eq:Da+1a}). Notice also that ${\mathcal R}_a \cap x^a L^2 \Omega^*$
is dense in ${\mathcal R}_a$. Therefore,
\[
\calD(D_{{\max},a}) = G_a \left({\mathcal R}_a \cap x^a L^2 \Omega^*\right)
+ {\mathcal N}_a, 
\]
where ${\mathcal N}_a$ is the nullspace of (\ref{eq:Daa-1}). 

Finally, if $\gamma \in \calH^*_{\min}(M,g,a)$, then $\langle D_{{\max},a} \eta, 
\gamma \rangle = 0$ for all $\eta \in \calD(D_{{\max},a})$, i.e.\ 
$\langle f, \gamma\rangle = 0$ for all $f \in {\mathcal R}_a \cap x^a L^2\Omega^*$.
Using the density statement above, this shows that $\gamma$ lies in the nullspace 
of (\ref{eq:Da+1a}), and is thus conormal.

We note in conclusion that any $\omega \in \calD(D_{{\max},a})$ has a `weak' 
asymptotic expansion of the form 
\[
\omega \sim \omega_0^-(y,z) x^{\gamma_0^-} + \omega_0^+(y,z) x^{\gamma_0^+} + \omega'
\]
where $\omega' = {\mathcal O}(x^{a+1})$, again in a suitable weak sense. We refer
to \ \cite{Ma-edge} for more details. 

\subsection{The minimal Hodge cohomology}
We are now in a position to prove the
\begin{theorem} Let $M$ be a manifold with an incomplete edge metric $g$. Then
\[
\calH^{k}_{{\min}}(M,g,a) =
\left\{ 
\begin{array}{rclll}
\mbox{Im}\,\left(I\!H^k_{\underline{\frakm} + <a>}(X,B) \right.
& \to &  \left. I\!H^k_{\overline{\frakm} + \ll a-1 \gg }(X,B)\right) & \quad & f 
\mbox{odd} \\
\mbox{Im}\,\left(I\!H^k_{\underline{\frakm} + <a-1/2>} (X,B) \right. & \to  & \left.
I\!H^k_{\overline{\frakm} + \ll a-1/2 \gg } (X,B)\right) & \quad & f \mbox{even}.
\end{array}
\right.
\]
In particular, when $a=0$, 
\[
\calH^{k}_{{\min}}(M,g,0) = \mbox{Im\,}(\IH^{k}_{\underline{\frakm}}(X,B) 
\longrightarrow \IH^{k}_{\overline{\frakm}}(X,B)).
\]
\end{theorem}
\begin{proof}
Recalling that $\calH^*_{{\min}}$ is quasi-isometry invariant, we may as well assume 
that $\Delta_F$ has no small nonzero eigenvalues, in the sense of the preceding subsections. 
If there are no small eigenvalues at all, i.e.\ either when $(a-(f+1)/2,a-(f-1)/2) \cap
{\mathbb N} = \emptyset$ (which holds, for example, when $a=0$ and $f$ is odd), 
or else if there exists $q_a \in (a-(f+1)/2,a-(f-1)/2)$ but $H^{q_a}(F) = \{0\}$, 
then the result follows directly from what we have already done, since then 
$H^k_{{\max}/{\min}}(M,g,a)$ and and $\calH^k_{{\max}/{\min}/\Kr/\Fr}(M,g,a)$ are all 
equal, cf. Corollary \ref{cor:trivext}. Thus we suppose that there exists 
$q_a \in (a-(f+1)/2,a-(f-1)/2)$ such that $H^{q_a}l(F) \neq \{0\}$. 

According to Theorem \ref{th:weightedco}, the space appearing on the right in the 
statement of this theorem, for $f$ even or odd, is identified with 
$\mbox{Im}\, \left(H^k_{{\min}}(M,g,a) \to H^k_{{\max}}(M,g,a)\right)$. For  
simplicity, we denote it as ${\mathcal J}^k(M,g,a)$.

We claim first that there is a natural injective map
\[
\calH^k_{{\min}}(M,g,a) \longrightarrow {\mathcal J}^k(M,g,a).
\]
To see this, recall that any form $\omega \in \calH^k_{{\min}}(M,g,a) = \calH^k_{\Kr}(M,g,a) 
\cap \calH^k_{\Fr}(M,g,a)$ naturally represents a class in $H^k_{{\min}}(M,g,a)$.
If $[\omega] = 0$ in $H^k_{{\max}}(M,g,a)$, then $\omega = d\zeta$ for some $\zeta \in 
\calD(d_{{\max},a})$. But this would imply that $||\omega||_a^2 = \langle \omega, 
d\zeta \rangle_a = 0$ since $\omega \in \ker(\delta_{{\min},a})$. 
This proves the claim.

The issue, then, is to prove that any class $[\eta] \in {\mathcal J}^k(M,g,a)$
is represented by an element of $\calH^k_{{\min}}(M,g,a)$. Choose a representative 
$\eta \in \calA^\gamma\Omega^k$, $\gamma > a - (f-1)/2$, for this class. Now use a 
generalized inverse $G$ for $D_a: x^a L^2\Omega^*(M,g) \to x^{a-1}L^2\Omega^*(M,g)$
(acting on its maximal domain). According to Propositions 8 and 9, this gives
$\zeta \in x^a L^2\Omega^* \cap \calA^{\gamma_0^-}$ and an element of the cokernel, which
by duality corresponds to an element $\omega \in \ker D_a \cap x^{a+1}L^2\Omega^* = 
\calH^k_{{\min}}(M,g,a)$, such that $\eta = D_a\zeta + \omega$. (Recall that $\gamma_0^-$ 
is the lower of the two indicial roots in the critical interval.) The theorem will follow 
once we show that $\delta_a \zeta = 0$. Now
\[
||\delta_a \zeta||_a^2 = \langle \delta_a \zeta, \eta - d\zeta - \omega\rangle_a.
\]
Integrating by parts formally this should vanish, so it remains to show that
each integration by parts is valid. First, $\langle \delta_a \zeta, \eta \rangle_a = 0$
since both terms are conormal and $d\eta = 0$. Similarly
$\langle \delta_a \zeta, \omega \rangle_a = 0$ since $\omega \in \ker d_{{\min},a}$.
To show that the remaining term vanishes, observe that $\zeta = 
x^{\gamma_0^-}\zeta_0^- + x^{\gamma_0^+} \zeta_0^+ + \zeta^{\prime}$ where
$\zeta^{\prime} \in \calA^{\gamma}\Omega^*$, $\gamma > a-(f-1)/2$.
Writing $\zeta_0^\pm = \mu_0^\pm + dx \wedge \nu_0^\pm$, then
$\mu_0^\pm, \nu_0^\pm \in \calC^\infty$ and in the nullspace of $\Delta_F$. 
A closer inspection of the equation $D_a \zeta =\eta - \omega$ shows that 
\[
\frakd \mu_0^\pm = \frakd \nu_0^\pm = \frakd^* \mu_0^\pm = \frakd^* \nu_0^\pm = 0.
\]
(The operators $\frakd$, $\frakd^*$ are the differential and codifferential
for the (fibre-harmonic) projected complex for the Riemannian submersion
metric $(Y,\tilde{g}(x))$.) In any case, the identity
$\langle d\zeta, \delta_a \zeta\rangle_a = 0$ is now immediate.
This proves the remaining assertion, and hence the theorem.

\end{proof}

\section{Hodge theory for complete edge metrics}
We are also able to determine the dimensions of the spaces of the spaces
$L^2\calH^k(M,g)$ when $(M,g)$ is a manifold with a complete edge metric. 
Unlike the incomplete case, in certain degrees this Hodge cohomology may
be infinite dimensional, i.e.\ there is an infinite dimensional space of 
$L^2$ harmonic forms. The simplest example of this is when $(M,g)$
is the $n$-dimensional hyperbolic space, or indeed any conformally
compact manifold, and $k = n/2$, cf.\ \cite{MaPh}. 

\begin{theorem}  Let $(M^n,g)$ be a manifold with a complete edge metric.
Let $X$ be the compact stratified space defined in \S 3. Suppose that $k$ is 
{\bf not} of the form $j + (b+1)/2$ where $\calH^j(F) \neq \{0\}$. Then
\[
L^2\calH^k(M,g) \cong {\IH}^k_{f+\frac{b}{2}-k}(X,B).
\]
In this case, the $L^2$ signature theorem for $M$ is the same as the 
$L^2$ signature theorem for $M$ endowed with 
the conformally equivalent incomplete edge metric $x^2 g$.
In all other cases, where $k$ does have this form, $L^2\calH^k(M,g)$ is 
infinite dimensional.
\end{theorem}

\begin{proof}
There are several viable ways to proceed: one could use a parametrix
construction based on the edge calculus to do a global Hodge theoretic
argument as in \cite{MaPh}; one could also, as in the incomplete case, use 
sheaf theory, calculations of local cohomologies, etc., ab initio; we take a 
shorter intermediate route, reducing to the incomplete edge case using the 
conformal invariance of the space of middle degree $L^2$ harmonic forms. 

It suffices to consider forms of degree $k \geq n/2$. Define $k = (n+r)/2$ with 
$0 \leq r \leq n$.  For later use, we also set $\sigma = (f+r)/2$, 
$\sigma' = (f-r)/2$. Now, since $2k = n + r$, $k$ is the middle degree on the manifold 
$\widetilde{M} = M \times S^{r}$. Endow $\widetilde{M}$ with the product metric 
$G_r = g + \kappa_{r}$ (where $\kappa_r$ is the standard metric on $S^{r}$); 
in a neighbourhood of $\del\widetilde{M}$
\[
G_r = \frac{dx^2 + h}{x^2} + (\kappa + \kappa_r).
\]
In other words, $G_r$ is still a complete edge metric with the same base $(B,h)$ as $g$, 
but with fibre $(F\times S^{r},\kappa + \kappa_{r})$. There is a corresponding 
{\it incomplete} edge metric $\hat{g} = x^2 \tilde{g}$ on $\widetilde{M}$. We 
denote its compactification, obtained by pinching the fibres $F \times S^{r}$ 
at the boundary, by $\widehat{X}$. 

We first claim that
\begin{equation}
L^2\calH^k(\widetilde{M},\tilde{g})=\calH^k_{\max}(\widetilde{M},\hat{g},0).
\label{eq:sss}
\end{equation}
The verification is straightforward; by conformal invariance of the $L^2$ condition
and the operator $\delta =\pm  *d*$ in the middle degree, harmonic
forms on the complete manifold are in the maximal domains of $d$ and $\delta$
(and indeed their nullspaces) on the incomplete manifold, and conversely.

By the $L^2$ K\"unneth theorem, 
\[
L^2\calH^k(\widetilde{M},\tilde{g}) = L^2\calH^k(M,g) \oplus L^2\calH^{k-r}(M,g).
\]
On the other hand, while we have not shown how to compute $\calH^k_{\max}
(\widetilde{M},\hat{g},0)$ in general, and indeed have noted that it is sometimes 
infinite dimensional, we now show that under certain hypotheses, it equals 
$\calH^k_{\min}(\widetilde{M},\hat{g},0)$; thus we can then apply Theorem 3 to
calculate the right side of (\ref{eq:sss}) as $\IH^k_{\frakm}(\widetilde{X})$. 
(Under these hypotheses, the intersection cohomologies with upper or lower middle
perversity are the same, so we just write $\frakm$.)

So, let us suppose that either $b$ is even, or else if $b$ is odd then $H^{\sigma}(F)$ 
(and hence $H^{\sigma'}(F)$) is trivial; note this last condition is automatic
when $r > f$. We claim that under these conditions, the minimal and maximal de Rham complexes 
on $\widetilde{M}$ coincide, and thus $\calH^k_{\max}(\widetilde{M},\hat{g},0) = 
\calH^k_{\min}(\widetilde{M},\hat{g},0)$. Indeed, this follows directly from 
Corollary~\ref{cor:trivext}: we have $\dim \widetilde{M} = \tilde{n} = 2k$ even; 
if $b$ is even, then $\tilde{f} = f + r = \tilde{n} - b  - 1$ is odd, while if 
$b$ is odd, then 
\[
H^{\tilde{f}/2}(F\times S^r) = H^{\sigma}(F) \oplus H^{\sigma'}(F) = \{0\}.
\]

Taking these facts together, and assuming this vanishing of the fibre cohomology
when $b$ is odd, we have proved that
\begin{equation}
L^2\calH^k(M,g) \oplus L^2\calH^{k-r}(M,g) = \IH^k_{\frakm}(\widetilde{X}).
\label{eq:ssss}
\end{equation}

It remains to compute the final term on the right of (\ref{eq:ssss}). We
decompose 
\[
\widetilde{X} = M\times S^r \sqcup_{Y \times S^r} \tilde{Z},
\]
where $\tilde{Z}$ is the cone bundle over $B$ with fibre $F \times S^r$ 
and boundary $Y \times S^r$. Of course,
\[
H^k(M\times S^r) = H^k(M) \oplus H^{k-r}(M).
\]
On the other hand, 
\[
\IH^*_{\frakm}(\tilde{Z}) = \IH^*_{\frakm}(Z) \oplus \IH^{*-r}_{\frakm}(Z).
\]
To prove this, we return to the sheaf-theoretic description. 
For product neighbourhoods $\calU = \calV \times C_1(F\times S^r)$, we have
\[
\IH^j_{\frakm}(\calU;\calL) = \IH^j_{\frakm}(C_1(F\times S^r)) = 
\left\{
\begin{array}{ccl} 
H^j(F) \oplus H^{j-r}(F),& \quad &j \leq (f+r-1)/2 \\
0 & &\mbox{otherwise.} 
\end{array} \right.
\]
Note that we have combined the conditional inequality in this last
step, which should depend on the parity of $f+r$, into one condition.
This condition is correct as stated if $f+r$ is odd; if $f+r$ is
even, then one would expect the condition $j \leq (f+r)/2 - 1 = \sigma-1$ or
$j \leq (f+r)/2 = \sigma$, depending on whether one was using upper or lower
middle perversity. However, the hypothesis $H^\sigma(F) = 0$ guarantees
that we get the same result in either case. Thus we see that the spectral
sequence whose hypercohomology computes the intersection cohomology of
$\tilde{Z}$ decouples into two noninteracting pieces.

We have now proved that for a fixed $k$, assuming the hypotheses above, 
\[
L^2\calH^k(M) \oplus L^2\calH^{k-r}(M) = \IH^k_{\frakm}(X) \oplus 
\IH^{k-r}_{\frakm}(X).
\]
We would, of course, like the summands to be equal separately; this can
be seen simply by noting that the correspondence takes place on the level of
forms, and we can separate out the terms with like degree.

It remains to show that in the remaining cases, i.e.\ when $b$ is odd
and $H^{\sigma}(F) \neq \{0\}$, $\sigma = (2k-b-1)/2$, $L^2\calH^k(M)$
is infinite dimensional. This follows from two assertions: that
$0$ is in the essential spectrum of $\Delta_k$, and that there is
a spectral gap at $0$ for $\Delta_k$. We content ourselves with sketching
the proofs briefly. The first step relies on the observation 
that up to quasiisometry, some neighbourhood of infinity
looks like the product of half of a hyperbolic space and a compact manifold,
specifically $\HH^{b+1}_+ \times F$. Since $b+1$ is even, there is 
an infinite dimensional family of $L^2$ harmonic forms on the first factor.
Since $k = (b+1)/2 + \sigma$, we can take suitable truncations of these, coupled with 
harmonic forms of degree $\sigma$ on $F$, to produce a Weyl sequence on $M$.
For the second step, we construct a parametrix in the edge calculus for 
$\Delta_k$. Its normal operator is given by
\[
N(\Delta_k) = \sum_{j=0}^k N(\Delta_{\HH^{b+1},j}) + \Delta_{F,k-j}.
\]
Both of these operators have spectral gaps at zero (even when $j=(b+1)/2$),
so it is possible to construct a parametrix for $\Delta_k - \lambda$
with compact remainder when $\lambda$ is small but nonzero. We refer
to \cite{Ma-conf.cpt.} and \cite{MaPh} for more complete descriptions
of such proofs in a slightly simpler context.
\end{proof}

\end{document}